\documentclass[10pt,twoside]{siamltex}
\usepackage{amsfonts,epsfig}

\usepackage{graphicx}
\usepackage{amssymb}
\usepackage{amsmath}

\usepackage{color}

\usepackage{eqparbox}

\usepackage{algorithm}
\usepackage{algorithmic}

\usepackage{xcolor}
\usepackage{tikz}

\definecolor{mygray}{RGB}{47,79,79}
\newcommand{\tc}{\textcolor{black}}

\def\e{{\varepsilon}}

\newcommand{\N}{\mathbb N}
\newcommand{\R}{\mathbb R}
\newcommand{\C}{\mathbb C}

\newcommand{\Sn}{\mathbb S}
\newcommand{\U}{\mathbb U}
\newcommand{\Z}{\mathbb Z}
\newcommand{\W}{\mathbb W}
\newcommand{\V}{\mathbb V}
\newcommand{\p}{r} 
\newcommand{\s}{s} 

\newcommand{\argmin}{\operatorname*{arg\,min}}
\newcommand{\colspan}{\operatorname{ran}}
\newcommand{\tr}{\operatorname{trace}}
\newcommand{\sym}{\operatorname{sym}}
\newcommand{\Skew}{\operatorname{skew}} 
\newcommand{\dist}{\operatorname{dist}}
\newcommand{\Exp}{\operatorname{Exp}}
\newcommand{\Log}{\operatorname{Log}}

\newcommand{\mcM}{\mathcal{M}}
\newcommand{\mcN}{\mathcal{N}}

\newtheorem{remark}{Remark}
\newtheorem{example}{Example}

\title{Manifold interpolation and model reduction}
\author{
  Ralf Zimmermann\thanks{Department of Mathematics and Computer Science, University of Southern Denmark (SDU) Odense,
    (zimmermann@imada.sdu.dk).}
}
\begin{document}
\maketitle

\begin{abstract}
    One approach to parametric and adaptive model reduction is via the interpolation of orthogonal bases, subspaces or positive definite system matrices.
    In all these cases, the sampled inputs stem from matrix sets that feature a geometric structure and thus form so-called matrix manifolds.
    This work will be featured as a chapter in the upcoming {\em Handbook on Model Order Reduction},
    (P. Benner, S. Grivet-Talocia, A. Quarteroni, G. Rozza, W.~H.~A. Schilders, L.~M. Silveira, eds, to appear on DE GRUYTER)
    and reviews the numerical treatment of the most important matrix manifolds that arise in the context of model reduction.
    Moreover, the principal approaches to data interpolation and Taylor-like extrapolation on matrix manifolds are outlined and complemented by algorithms in pseudo-code.
\end{abstract}


\section{Introduction \& Motivation}
\label{sec:intro}
This work addresses interpolation approaches for parametric model reduction.
This includes techniques for 
\begin{itemize}
 \item computing trajectories of parameterized subspaces,
 \item computing trajectories of parameterized reduced orthogonal bases,
 \item structure-preserving interpolation.
\end{itemize}
Mathematically, this requires data processing on nonlinear matrix manifolds.
The exposition at hand intends to be an introduction and a reference guide 
to numerical procedures with matrix manifold-valued data.
As such it addresses practitioners and scientists new to the field.
It covers the essentials of those matrix manifolds that arise most frequently in practical problems in model reduction.
The main purpose is not to discuss concrete model reduction applications, but rather to provide the essential tools, building blocks and background theory to enable the reader to devise her/his own approaches for such applications.

The text was designed such that it works as a commented formula collection, meanwhile giving sufficient context, explanations and, not least, precise references to enable the interested reader to  immerse further in the topic.

\subsection{Parametric model reduction via manifold interpolation: An introductory example}
\label{sec:introInterp4pMOR}
The basic objective in model reduction is to emulate a large-scale dynamical system
with very few degrees of freedom
such that its input/output behavior is preserved as well as possible.
While classical model reduction techniques aim at producing an accurate 
low-order approximation to the \tc{autonomous behavior} of the original system,
parametric model reduction (pMOR) tries to account for additional system parameters.
\tc{If we look for instance at aircraft aerodynamics, an important task is to solve the unsteady Navier-Stokes equations at various flight conditions, which are, amongst others, specified by the altitude, the viscosity of the fluid (i.e. the Reynolds number) and the 
relative velocity (i.e. the Mach number).}
We explain the objective of pMOR
with the aid of a generic example 
\tc{ in the context of proper orthogonal decomposition-based model reduction.
Similar considerations apply to frequency domain approaches, Krylov subspace methods and balanced truncation, which are discussed in other chapters of the upcoming Handbook on Model Order Reduction.
}
Consider a spatio-temporal dynamical system in semi-discrete form
\begin{equation}
\label{eq:basicPDE}
  \frac{\partial}{\partial t}x(t, \mu) = f(x(t, \mu); \mu), \quad x(t_0,\mu) = x_{0,\mu},
\end{equation}
where $x(t, \mu)\in\R^n$ is the spatially discretized {\em state vector} of dimension $n$,
the vector $\mu=(\mu_1,\ldots,\mu_d)\in \R^d$ accounts for additional system parameters 
and $f(\hspace{0.1cm} \cdot\hspace{0.1cm}; \mu):\R^n \rightarrow \R^n$ is the (possibly nonlinear, parameter-dependent) right hand side function.
Projection-based MOR starts with constructing a suitable low-dimensional subspace
that acts as a space of candidate solutions.

{\em Subspace construction.}
One way to construct the required projection subspace is the
proper orthogonal decomposition (POD), \cite{HinzeVolkwein2005}.
In its simplest form, the POD can be summarized as follows.
For a fixed system parameter $\mu=\mu_0$, let
$x^1:=x(t_1,\mu_0),...,x^{m}:=x(t_m,\mu_0) \in \R^n$ 
be a set of state vectors satisfying \eqref{eq:basicPDE} and let $\Sn:= \left(x^1,...,x^{m}\right)\in \R^{n\times m}$.
The state vectors $x^i$ are called 
{\em snapshots} and the matrix $\Sn$ is called the associated {\em snapshot matrix}.
POD is concerned with finding a subspace $\mathcal{V}$ of dimension $\p \leq m$
represented by a column-orthogonal matrix $\V_r\in \R^{n\times \p}$
such that the error between the input snapshots
and their orthogonal projection onto $\mathcal{V}=\colspan(\V_r)$ is minimized:
\[
  \min_{V\in \R^{n\times \p}, V^TV=I}{\sum_k{\|x^k -VV^Tx^k \|^2_2}} \quad \left(\Leftrightarrow   
  \min_{V\in \R^{n\times \p}, V^TV=I}{\|\Sn - VV^T\Sn\|^2_F } \right).
\]
The main result of POD is that for any $\p\leq m$, 
the best $\p$-dimensional approximation of $\colspan(x^1,...,x^{m})$ in the above sense
is $\mathcal{V}=\colspan(v^1,...,v^{\p})$,
where $\{v^1,...,v^{\p}\}$ are the eigenvectors of the matrix $\Sn\Sn^T$
corresponding to the $\p$ largest eigenvalues.
The subspace $\mathcal{V}$ is called the {\em POD subspace}
and the matrix $\V_r = (v^1,...,v^{\p})$ is the {\em POD basis matrix}.
The same subspace is obtained via a compact singular value decomposition (SVD) 
of the snapshot matrix $\Sn = \V \Sigma \Z^T$, truncated to the first
$\p\leq m$ columns of $\V\in\R^{n\times m}$ by setting $\mathcal{V}:=\colspan(\V_{\p})$. For more details, see, e.g. \cite[\S 3.3]{BennerGugercinWillcox2015}.
In the following, we drop the index $\p$ and assume that $\V$ is already the truncated matrix $\V = (v^1,...,v^{\p})\in\R^{n\times \p}$.

Since the input snapshots are supplied at a fixed system parameter vector $\mu_0$,
the POD subspace is considered to be an appropriate
space of solution candidates $\mathcal{V}(\mu_0)=\colspan(\V(\mu_0))$
at $\mu_0$.

{\em Projection.}
POD leads to a parameter decoupling
\begin{equation}
 \label{eq:PODansatz}
  \tilde{x}(t, \mu_0)= \V(\mu_0)x_{\p}(t).
\end{equation}
In this way, the time trajectory of the reduced model is uniquely defined
by the coefficient vector $x_{\p}(t)\in \R^\p$
that represents the reduced state vector with respect to the
subspace $\colspan(\V(\mu_0))$.
Given a matrix $\W(\mu_0)$ such that the matrix pair $\V(\mu_0), \W(\mu_0)$
is bi-orthogonal, i.e. $\W(\mu_0)^T\V(\mu_0) = I$,
the original system \eqref{eq:basicPDE} can be reduced in dimension as follows.
Substituting \eqref{eq:PODansatz} in \eqref{eq:basicPDE}
and multiplying with $\W(\mu_0)^T$ from the left leads to
\begin{equation}
\label{eq:basicMOR}
       \frac{d}{dt}x_{\p}(t) = \W^T(\mu_0) f(\V(\mu_0) x_{\p}(t); \mu_0), \quad x_\p(t_0) = \V^T(\mu_0)x_{0, \mu_0}.
\end{equation}
This approach goes by the name of Petrov-Galerkin projection,
if $\W(\mu_0)\neq \V(\mu_0)$ and Galerkin projection if $\W(\mu_0)=\V(\mu_0)$.
There are various ways to proceed from \eqref{eq:basicMOR}
depending on the nature of the function $f$ and many of them are discussed in other chapters of the upcoming
Handbook on Model Order Reduction.
\footnote{
If $f(\hspace{0.1cm} \cdot\hspace{0.1cm}; \mu_0)$ is linear, the reduced operator 
$\W^T(\mu_0) \circ f(\hspace{0.1cm} \cdot\hspace{0.1cm}; \mu_0) \circ \V(\mu_0)$
can be computed a priori ({\em `offline'}) and stays fixed throughout the
time integration.
If $f(\hspace{0.1cm} \cdot\hspace{0.1cm}; \mu_0)$ is affine, the same approach can be carried over to the affine
building blocks of $f(\hspace{0.1cm} \cdot\hspace{0.1cm}; \mu_0)$, see e.g. \cite{HaasdonkOhlberger2011}.
For a nonlinear $f(\hspace{0.1cm} \cdot\hspace{0.1cm}; \mu_0)$, an affine approximation can be constructed via the
emperical interpolation method (EIM, \cite{Barrault2004}).
Other approaches that address nonlinearities include the
discrete empirical interpolation method (DEIM, \cite{DEIM_Chaturantabut_2010})
and the missing point estimation (MPE, \cite{MPE_Astrid2008, ZimmermannWillcox2016}).}

For illustration purposes, we proceed with $\W(\mu_0)=\V(\mu_0)$
and assume that the right hand side function $f$ splits into a linear and a nonlinear part:
$f(x;\mu_0) = A(\mu_0)x + \mathbf{f}(x;\mu_0)$, where $A(\mu_0) \in\R^{n\times n}$ is, say, a \tc{symmetric and negative definite matrix to foster stability.}
Then, \eqref{eq:basicMOR} becomes
\[
  \frac{d}{dt}x_{\p}(t) 
  = \V^T(\mu_0) A(\mu_0) \V(\mu_0)x_{\p}(t) + \V^T(\mu_0)\mathbf{f}\bigl(\V(\mu_0) x_{\p}(t); \mu_0\bigr).
\]
In the discrete empirical interpolation method (DEIM, \cite{DEIM_Chaturantabut_2010}), the large-scale nonlinear term $\mathbf{f}\bigl(\V(\mu_0)x_{\p}(t);\mu_0)$ is approximated via a mask matrix $P = (e_{i_1},\ldots,e_{i_\s})\in \R^{n\times \s}$,
\tc{where $\{i_1,\ldots,i_\s\}\subset \{1,\ldots,n\}$ and $e_j = (\ldots,\stackrel{j}{1},\ldots)^T\in \R^n$ is the $j$th canonical unit vector.}
The mask matrix $P$ acts as an entry selector on a given $n$-vector via $P^T v = (v_{i_1},\ldots,v_{i_\s})^T\in \R^\s$.
In addition, another POD basis matrix $\U(\mu_0)\in\R^{n\times \s}$ is used, which is obtained from snapshots of the nonlinear term.
The matrices $P$ and $\U(\mu_0)$ are combined to form an oblique projection of the non-linear term onto the subspace $\colspan(\U(\mu_0))$.
This leads to the reduced model 
\begin{eqnarray}
\nonumber
       \frac{d}{dt}x_{\p}(t)& =&
       \V^T(\mu_0) A(\mu_0) \V(\mu_0)x_{\p}(t)\\
       \label{eq:basicDEIM}
       &&
       + \V^T(\mu_0)\U(\mu_0)(P^T\U(\mu_0))^{-1}P^T\mathbf{f}\bigl(\V(\mu_0) x_{\p}(t); \mu_0\bigr),
\end{eqnarray}
whose computational complexity is formally independent of the full-order dimension $n$, see \cite{DEIM_Chaturantabut_2010} for details.
Mind that by assumption, $M(\mu_0):=-\V^T(\mu_0) A(\mu_0) \V(\mu_0)$ is symmetric positive definite and that
both $\V(\mu_0)$ and $\U(\mu_0)$ are column-orthogonal.
\tc{Moreover, for a fixed mask matrix $P$, coordinate changes of $\V(\mu_0)$ and $\U(\mu_0)$ do not affect the approximated state
$\tilde{x}(t, \mu_0)= \V(\mu_0)x_{\p}(t)$, so that essentially, the reduced system
\eqref{eq:basicDEIM} depends only on the subspaces $\colspan(\V(\mu_0))$ and $\colspan(\U(\mu_0))$
rather than the matrices $\V(\mu_0)$ and $\U(\mu_0)$.}\footnote{\tc{Replacing $\U$ with $\U S$, $S\in\R^{\s\times\s}$ orthogonal, does not affect \eqref{eq:basicDEIM} at all.
Replacing $\V$ with $\V R$, $R\in\R^{\p\times\p}$ orthogonal, induces a coordinate change on the reduced state $ x_\p = R \hat{x}_\p$ but preserves the output
$\tilde{x}(t)= \V x_{\p}(t) = \V R \hat{x}_\p(t)$.}}

Solving \eqref{eq:basicMOR}, \eqref{eq:basicDEIM} constitutes the {\em online stage} of model reduction.
The main focus of this exposition is {\em not} on the efficient solution of the 
reduced systems \eqref{eq:basicMOR} or \eqref{eq:basicDEIM}
at a fixed $\mu_0$, but on
tackling parametric variations in $\mu$.
In view of the associated computational costs, 
it is important that this can be achieved {\em without} computing additional snapshots in the online stage.

\tc{A straightforward way to achieve this is to extend the snapshot sampling to the $\mu$-parameter range
to produce POD basis matrices that are to cover all input parameters.
This is usually referred to as the ``global approach''. For nonlinear systems, the global approach may suffer from requiring a large number of snapshot samples. Moreover, the snapshot information is blurred in the global POD and features that occur only in a restricted regime affect the ROM predictions everywhere. Therefore, localized approaches are preferable,
see e.g. 
\cite{FranzZimmermannGoertzKarcher2014, OS15, LDEIM, SargsyanBK18, Zimmermann2014}.}

In this contribution, the focus is on constructing trajectories of 
functions in the system parameters $\mu$ on certain sets of structured matrix spaces.
In the above example, these are the symmetric positive definite matrices $\{M\in\R^{\p\times \p}| M^T = M, v^TMv>0\hspace{0.1cm}  \forall v\neq 0\}$, 
the orthonormal basis matrices $\{U\in \R^{n\times \s}| U^TU=I\}$ or the associated $\s$-dimensional subspaces 
$\mathcal{U}:= \colspan(U)\subset \R^n$:
\begin{eqnarray*}
\label{eq:spd_traj}
\mu \mapsto -\V^T(\mu) A(\mu) \V(\mu)&\in& \{M\in\R^{\p\times \p}| M^T = M, v^TMv>0\hspace{0.1cm} \forall v\neq 0\},\\
\label{eq:subspace_traj}
 \mu \mapsto \U(\mu)&\in& \{U\in \R^{n\times \s}| U^TU=I\},
 \\ 
 \mu \mapsto \mathcal{U}(\mu) = \colspan(\U(\mu)) &\in& \{\mathcal{U}\subset \R^n| \hspace{0.1cm}\mathcal{U} \text{ subspace}, \dim(\mathcal{U})= \s\}.
\end{eqnarray*}
%
We outline generic methods for constructing such trajectories via interpolation.
All the special sets of matrices considered above feature a differentiable structure that allows to consider them as submanifolds of some Euclidean matrix space, referred to as matrix manifolds.
The above example is not exhaustive. Other matrix manifolds may arise in model reduction applications.
To keep the exposition both general and modular, the interpolation techniques will be formulated for arbitrary submanifolds.
Model reduction literature on manifold interpolation problems includes
\cite{AmsallemFarhat2008, AmsallemFarhat2011, BennerGugercinWillcox2015, 
Degroote_etal2010, Morzi2007, Nguyen2012, Tadmor2011, Panzer_etal2010, Zimmermann2014, Choi2015, Massart2018}.

\subsection{Structure and organization}
The text is constructed modular rather than consecutive, so that selected reading is enabled. Yet, this entails that the reader will encounter some repetition.\\
Section \ref{sec:diffgeo_intro} covers the essential background from differential geometry.
Section \ref{sec:ManInterp} contains generic methods for interpolation and extrapolation on matrix manifolds.
In Section \ref{sec:MatMan}, the geometric and numerical aspects of the matrix manifolds that arise most frequently in the context of model reduction are discussed.\\
A practitioner that faces a problem in matrix manifold interpolation may skim through the recap on elementary differential geometry in Section \ref{sec:diffgeo_intro} and then move on to the appropriate subsection of Section \ref{sec:MatMan} that corresponds to the matrix manifold in the application. This provides the specific ingredients and formulas for conducting the generic interpolation methods of Section \ref{sec:ManInterp}.

\subsection{Notation \& Abbreviations} 
%
\begin{small}
\begin{itemize}
\item w.r.t.: with respect to
\item 
 EVD: eigenvalue decomposition
\item
 SVD: singular value decomposition
\item POD: proper orthogonal decomposition
 \item
 LTI: linear time-invariant (system)
 \item
 ODE: ordinary differential equation
 \item
 PDE: partial differential equation
 \item ONB: orthonormal basis
 \item 
 $\R^{n\times \p}$: the set of real $n$-by-$\p$ matrices
 \item
 $I_n$: the $n$-by-$n$ identity matrix; if dimensions are clear, written as $I$
 \item
 $\colspan(A)$: the subspace spanned by the columns of $A\in \R^{n\times \p}$
\item
 $GL(n)$: the general linear group of real, invertible $n$-by-$n$ matrices
\item
 $\sym(n)= \{A\in\R^{n\times n}| A^T =A\}$: the set of real, symmetric $n$-by-$n$ matrices
\item
$\Skew(n)= \{A\in\R^{n\times n}| A^T =-A\}$: the set of real, skew-symmetric $n$-by-$n$ matrices
\item
 $SPD(n)= \{A\in\sym(n)| x^TAx> 0 \forall x\in\R^n\setminus\{0\}\}$: the set of real, symmetric positive definite  $n$-by-$n$ matrices
 \item
 $O(n)= \{Q\in\R^{n\times n}| Q^TQ=I_n=QQ^T\}$: the orthogonal group
 \item
  $SO(n)= \{Q\in O(n)| \det(Q) = 1\}$: the special orthogonal group 
\item
 $St(n,\p)= \{U\in\R^{n\times \p}| U^TU=I_\p\}$: the (compact) Stiefel manifold, $\p\leq n$
\item
 $Gr(n,\p)$: the Grassmann manifold of $\p$-dimensional subspaces of $\R^n$, $\p\leq n$
 \item
 $\mcM$: a differentiable manifold
 \item $\mathcal{D}_p\subset \mcM$: an open domain around the point $p$ on a manifold $\mcM$
 \item $D_x\subset \R^n$: an open domain in the Euclidean space around a point $x\in\R^n$
\item
$T_p\mcM$: the tangent space of $\mcM$ at a location $p\in\mcM$
\item
 $\langle A,B\rangle_0 = \tr(A^TB)$: the standard (Frobenius) inner product on $\R^{n\times \p}$
\item
 $\langle v,w\rangle_p^\mcM$: the Riemannian metric on $T_p\mcM$ (the superscript is often omitted)
\item
 $\exp_m$: standard matrix exponential
 \item
 $\log_m$: standard (principal) matrix logarithm
\item 
$\Exp_p^\mcM$: the Riemmanian exponential of a manifold $\mcM$ at base point $p\in\mcM$
\item
$\Log_p^\mcM$: the Riemmanian logarithm of a manifold $\mcM$ at base point $p\in\mcM$
\end{itemize}
\end{small}
%
%
%
\section{Basic concepts of differential geometry}
\label{sec:diffgeo_intro}
This section provides the essentials on elementary differential geometry.
Established textbook references on differential geometry include \cite{DoCarmo2013riemannian, KobayashiNomizu1963, KobayashiNomizu1969, kuhnel2006differential, Lee1997riemannian}; condensed introductions can be found in \cite[Appendices C.3, C.4, C.5]{HelmkeMoore1994} and \cite{Gallier2011}. An account of differential geometry that is tailor-made to matrix manifold applications is given in \cite{AbsilMahonySepulchre2008}.

The fundamental objects of study in differential geometry are {\em differentiable manifolds}.
Differentiable manifolds are generalizations of curves (one-dimensional) and surfaces (two-dimensional)
to arbitrary dimensions.
Loosely speaking, an $n$-dimensional differentiable manifold $\mcM$ is a topological space that `locally looks like $\R^n$'
with certain smoothness properties.
This concept is rendered precisely by postulating that for every point $p\in \mcM$,
there exists a so-called {\em coordinate chart} $x:\mcM\supset \mathcal{D}_p \rightarrow \R^n$ that bijectively maps an open neighborhood $\mathcal{D}_p\subset \mcM$ of a location $p$ to an open neighborhood $D_{x(p)}\subset \R^n$ around $x(p)\in \R^n$ 
with the important additional property that the {\em coordinate change} 
$$x\circ \tilde{x}^{-1}: \tilde{x}(\mathcal{D}_p\cap\tilde{\mathcal{D}}_p) \rightarrow x(\mathcal{D}_p\cap\tilde{\mathcal{D}}_p)$$
of two such charts
$x,\tilde{x}$ is a diffeomorphism, where their domains of definition overlap, see \cite[Fig. 18.2, p. 496]{Gallier2011}
or \cite[Fig. 3.1, p. 342]{HelmkeMoore1994}.
Note that the coordinate change $x\circ \tilde{x}^{-1}$
maps from an open domain of $\R^n$ to an  open domain of $\R^n$, so that the standard concepts of multivariate calculus apply.
For details, see \cite[\S 3.1.1]{AbsilMahonySepulchre2008} or \cite[\S 18.8]{Gallier2011}.
Depending on the context, we will write $x(p)$ for the value of a coordinate chart at $p$ and also
$x\in\R^n$ for a point in $\R^n$.

Of special importance to numerical applications are {\em embedded submanifolds in the Euclidean space}.
\begin{definition}[Submanifolds of $\R^{n+d}$]
 A {\em parameterization} is an bijective differentiable function $f: \R^n \supset D \rightarrow f(D)\subset \R^{n+d}$ with continuous inverse such that its Jacobi matrix $Df_x\in \R^{(n+d)\times n}$ has full rank $n$ at every point $x\in D$.
 
 A subset $\mcM \subset \R^{n+d}$ is called an {\em $n$-dimensional embedded submanifold of $\R^{n+d}$},
 if for every $p\in \mcM$, there exists an open neighborhood $\Omega\subset \R^{n+d}$ such that 
 $\mathcal{D}_p:=\mcM\cap \Omega$ is the image of a parameterization 
 \[
  f: \R^n\supset D_x \rightarrow f(D_x) = \mathcal{D}_p = \mcM\cap \Omega \subset \R^{n+d}.
 \]
\end{definition}
One can show that if $f: D\rightarrow \mcM\cap \Omega$ and $\tilde{f}: \tilde{D}\rightarrow \mcM\cap \tilde{\Omega}$ are two parameterizations, say with $f(x_0)  = \tilde{f}(\tilde{x}_0)= p\in \mcM\cap \Omega\cap \tilde{\Omega}$, then 
\[
 \left(f^{-1}\circ \tilde{f}\right): \tilde{f}^{-1}(\Omega\cap \tilde{\Omega}) \rightarrow f^{-1}(\Omega\cap \tilde{\Omega})
\]
is a diffeomorphism (between open sets in $\R^n$). In this sense, parameterizations $f$ are the inverses of coordinate charts $x$.
In addition to coordinate charts and parameterizations, submanifolds can be characterized via \tc{equality constraints}.
This fact is due to the inverse function theorem of classical multivariate calculus \cite[\S I.5]{lang2001fundamentals}.
For details, see \cite[Thm. 18.7, p. 497]{Gallier2011}.
\begin{theorem}[{\cite[Prop. 18.7, p. 500]{Gallier2011}}]
\label{thm:regurbild}
 Let $h: \R^{n+d}\supset \Omega \rightarrow \R^{d}$ be differentiable and $c_0\in \R^d$ be defined such that 
 the differential $Dh_p\in \R^{d\times (n+d)}$ has maximum possible rank $d$ at every point $p\in \Omega$ with $h(p) = c_0$.
 Then, the preimage
 \[
  h^{-1}(c_0) = \{p\in \Omega| \hspace{0.1cm} h(p) = c_0\}
 \]
 is an $n$-dimensional submanifold of $\R^{n+d}$.
\end{theorem}
An obvious application of Theorem \ref{thm:regurbild} to the function $h: \R^3 \rightarrow \R, (x_1,x_2,x_3) \mapsto x_1^2+x_2^2 + x_3^2 -1$
establishes the  unit sphere $S^2 = h^{-1}(0)$ as a $2$-dimensional submanifold of $\R^{2+1}$.
As a more sophisticated example, we recognize the orthogonal group as a differentiable (sub)-manifold:
\begin{example}
\label{ex:reg_urbild_On}
 Consider the orthogonal group $O(n)\subset \R^{n\times n}\simeq \R^{n^2}$
 and the set of symmetric matrices $\sym(n)\simeq \R^{n(n+1)/2}$.
 Define $h: \R^{n\times n}\rightarrow \sym(n), A \mapsto A^TA -I$.
 Then $Dh_A(B) = A^TB + B^T A$. For $Q\in O(n)$, the differential is indeed surjective: For any $M\in \sym(n)$, it holds
 $Dh_Q( \frac12 QM) = \frac12 Q^TQM + \frac12 M^TQ^TQ = M$.
 As a consequence, the orthogonal group $O(n)$ is a submanifold of dimension 
 $n^2 - \frac12(n(n+1)) = \frac12(n(n-1))$ of the Euclidean matrix space $\R^{n\times n}$.
\end{example}

\subsection{Intrinsic and extrinsic coordinates.}
\label{sec:extr_intr_coordinates}
As a rule, numerical data processing on manifolds requires calculations in explicit coordinates. 
For differentiable submanifolds, we distinguish between two types: {\em extrinsic} and {\em intrinsic coordinates}.
Extrinsic coordinates address points on a submanifold  $\mcM\subseteq \R^n$ with respect to their coordinates in the ambient space $\R^n$, while intrinsic coordinates are with respect to the local parameterizations.
Hence, extrinsic coordinates are what an outside observer would see, while intrinsic coordinates correspond to the perspective of an observer that resides on the manifold.
Let's exemplify these concepts on the two-dimensional unit sphere $S^2$, embedded in $\R^3$.
As a point set, the sphere is defined by the equation
\[
 S^2 = \{(x_1,x_2,x_3)^T\in \R^3| \hspace{0.1cm}  x_1^2+x_2^2+x_3^2 = 1\}.
\]
Any three-vector $(x_1,x_2,x_3)^T\in S^2$ specifies a point on the sphere in {\em extrinsic} coordinates.
However, it is intuitively clear that $S^2$ is intrinsically a two-dimensional object.
Indeed, $S^2$  can be parameterized via
\[
 f:  \R^2 \supset [0,2\pi)^2 \rightarrow S^2\subset  \R^3, \quad (\alpha,\beta) \mapsto 
 \begin{pmatrix}
  \sin(\alpha)\cos(\beta)\\
  \sin(\alpha)\sin(\beta)\\
  \cos(\alpha)
 \end{pmatrix}.
\]
The parameter vector $(\alpha,\beta)\in \R^2$ specifies a point on $S^2$ in {\em intrinsic} coordinates.
Even though intrinsic coordinates directly reflect the dimension of the manifold at hand,
they often cannot be calculated explicitly and
extrinsic coordinates are the preferred choice in numerical applications \cite[\S 2, p. 305]{EdelmanAriasSmith1999}.
Turning back to Example \ref{ex:reg_urbild_On}, we recall that the intrinsic dimension of the orthogonal group is
$\frac12n(n-1)$. Yet, in practice, one uses the extrinsic representation with $(n\times n)$-matrices $Q$,
keeping the defining equation $Q^TQ=I$ in mind.

\subsection{Tangent spaces.}
\label{sec:tangspaces}
We need a few more fundamental concepts.
\begin{definition}[Tangent space  of a differentiable submanifold]
\label{def:tangentspace}
 Let $\mcM \subset \R^{n+d}$ be an $n$-dimensional submanifold of $\R^{n+d}$.
 The {\em tangent space} of $\mcM$ at a point $p\in \mcM$, in symbols $T_p\mcM$, is the space of 
 velocity vectors of differentiable curves $c:t \mapsto c(t)$ passing through $p$, i.e.,
 \[ 
  T_p\mcM = \{\dot{c}(t_0)| \hspace{0.1cm}  c:J\rightarrow \mcM,   \hspace{0.1cm} c(t_0)=p\}.
 \]
 \tc{Here, $J\subseteq \R$ is an arbitrarily small open interval with $t_0\in J$.}
\end{definition}
\begin{figure}[ht]
\centering
\includegraphics[width=1.0\textwidth]{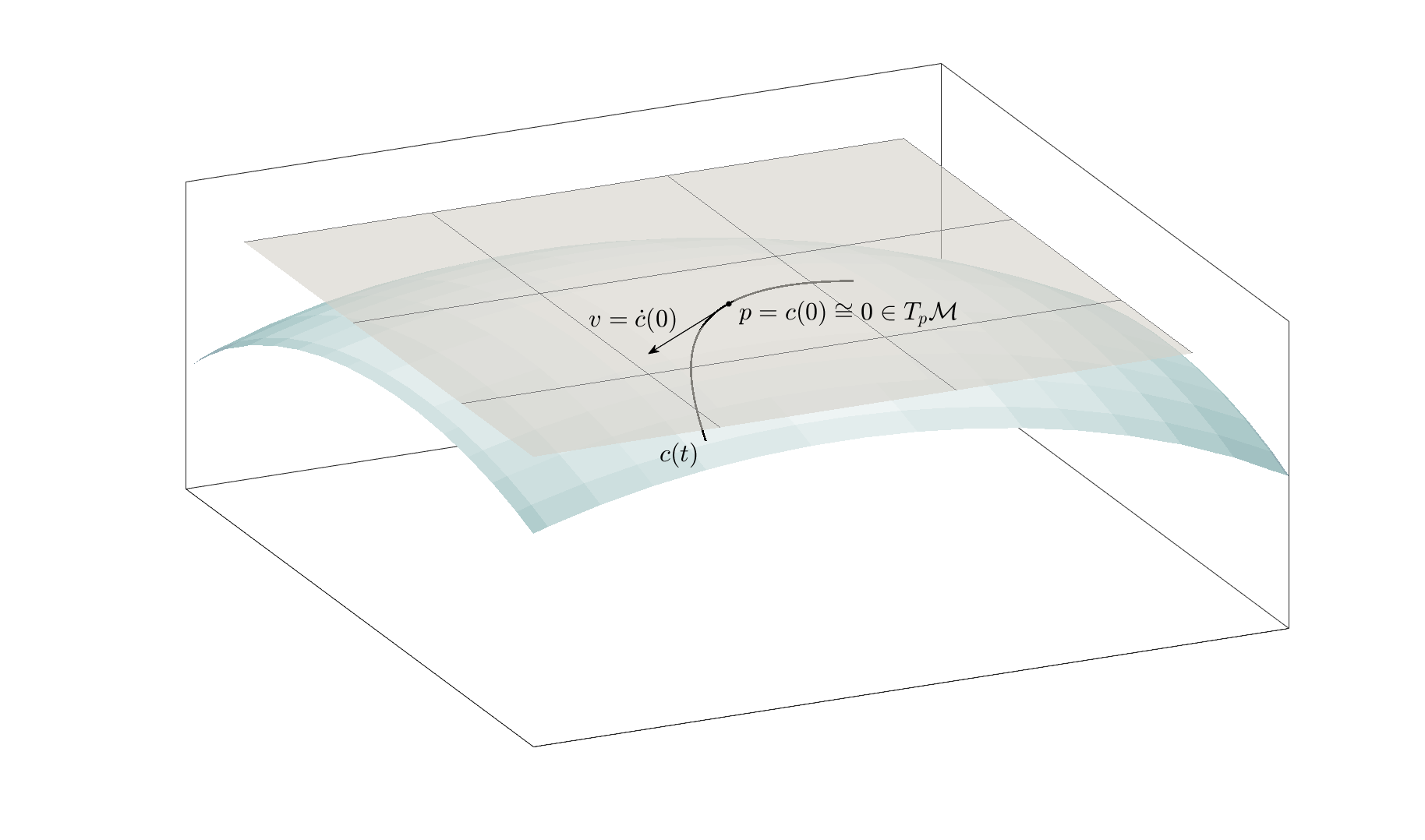}
\caption{Visualization of a manifold (curved surface)
with the tangent space $T_p\mcM$ attached.
The tangent vector $v = \dot{c}(0)\in T_p\mcM$ is the velocity vector of a curve $c:t\mapsto c(t)\in \mcM$.
}
\label{fig:manifold_plot_diffgeo}
\end{figure}
%
It is straightforward to show that the tangent space is actually a vector space.
Moreover, the tangent space can be characterized both with respect to intrinsic and extrinsic coordinates.
\begin{theorem}[Tangent space, intrinsic characterization]
 \label{thm:tangentspace_intr}
  Let $\mcM \subset \R^{n+d}$ be an $n$-dimensional submanifold of $\R^{n+d}$ and let
  $f:\R^n \supseteq D \rightarrow f(D)\subseteq \mcM$ be a parameterization.
  Then, for $x\in D$ with $p = f(x)\in\mcM$, it holds
  \[
   T_p\mcM = \colspan(Df_x).
  \]
\end{theorem}
\begin{theorem}[Tangent space, extrinsic characterization]
 \label{thm:tangentspace_extr}
  Let $h: \R^{n+d}\supset \Omega \rightarrow \R^{d}$  and  $c_0\in \R^d$ be as in Theorem \ref{thm:regurbild}
  and let $\mcM := h^{-1}(c_0) \subset \R^{n+d}$.
  Then, for  $p \in\mcM$, it holds
  \[
   T_p\mcM = \ker(Dh_p).
  \]
\end{theorem}
Note that both Theorem \ref{thm:tangentspace_intr} and Theorem \ref{thm:tangentspace_extr} immediately show that the tangent space
$T_p\mcM$ is a vector space of the same dimension $n$ as the manifold $\mcM$.
\begin{example}
 \label{ex:tangspace_On}
 The tangent space of the orthogonal group $O(n)$ at a point $Q_0$ is
 \[
  T_{Q_0}O(n) = \{\Delta \in \R^{n\times n}|\hspace{0.2cm} \Delta^TQ_0 = -Q_0^T \Delta\}.
 \]
 This fact can be established via considering a matrix curve $Q:t\mapsto Q(t)$ with $Q(0) = Q_0$
 and velocity vector $\Delta = \dot{Q}(0)\in T_{Q_0}O(n)$.
 Then, 
 \[
  0 =\frac{d}{dt}|_{t=0}I =   \frac{d}{dt}|_{t=0} Q^T(t)Q(t) = \Delta^TQ_0 + Q_0^T\Delta.
 \]
 (The claim follows by counting the dimension of the subspace $\{\Delta^TQ_0 = -Q_0^T \Delta\}$.)
 As an alternative, we can consider $h: \R^{n\times n}\rightarrow \sym, A \mapsto A^TA -I$ as in Example \ref{ex:reg_urbild_On}.
 Then $Dh_{Q_0}(\Delta) = Q_0^T\Delta + \Delta^T Q_0$ and $T_{Q_0}O(n) = \ker(Dh_{Q_0})$.
\end{example}
\subsection{Geodesics and the Riemannian distance function}
\label{sec:diffgeo_geodesics}
One of the most important problems in both general differential geometry 
and data processing on manifolds is to determine the shortest connection between two points on a given manifold. This requires to measure the lengths of curves.
Recall that the length of a curve $c:[a,b] \rightarrow \R^n$ in the Euclidean space is $L(c) = \int_a^b \|\dot{c}(t)\| dt$.
In order to transfer this to the manifold setting, an inner product for tangent vectors is needed that is consistent with the manifold structure.
\begin{definition}[Riemannian metrics]
  \label{def:RiemannMetric} 
  Let $\mcM$ be a differentiable submanifold of $\R^{n+d}$.
  A {\em Riemannian metric} on $\mcM$ is a family $(\langle\cdot,\cdot\rangle_p)_{p\in\mcM}$ of inner products $\langle\cdot,\cdot \rangle_p: T_p\mcM\times T_p\mcM\rightarrow \R$ that is smooth in variations of the base point $p$.\\
  The {\em length} of a tangent vector $v\in T_p\mcM$ is $\|v\|_p := \sqrt{\langle v,v\rangle_p}$.\footnote{This notation should not be confused with the classical $p$-norm $\sqrt[p]{\sum_i |v_i|^p}$.}
  The length of a curve $c:[a,b] \rightarrow \mcM$ is defined as $$L(c) = \int_a^b \|\dot{c}(t)\|_{c(t)} dt  = \int_a^b \sqrt{\langle\dot{c}(t),\dot{c}(t)\rangle_{c(t)}} dt.$$\\
  A curve is said to be {\em parameterized by the arc length}, if $L(c|_{[a,t]}) = t-a$ for all $t\in [a,b]$. Obviously, {\em unit-speed curves} with $\|\dot{c}(t)\|_{c(t)}\equiv 1$ are parameterized by the arc length.
  Constant-speed curves with $\|\dot{c}(t)\|_{c(t)}\equiv \nu_0$ are parameterized proportional to the arc length.
  The {\em Riemannian distance} between two points $p,q\in \mcM$ with respect to a given metric is
  \begin{equation}
\label{eq:RiemannDist}
  \dist_{\mcM}(p,q) = \inf\{L(c)| c:[a,b]\rightarrow \mcM \mbox{ piecewise smooth, } c(a)=p, c(b)=q\},
  \end{equation}
  where, by convention, $\inf\{\emptyset\} =\infty$. 
\end{definition}
Hence, a shortest path between $p,q\in \mcM$ is a curve $c$ that connects $p$ and $q$ such that 
$L(c) = \dist_{\mcM}(p,q)$.
In general, shortest paths on $\mcM$ do not exist.\footnote{Consider $\R^{2,*} = \R^2\setminus\{(0,0)\}$ with the Euclidean inner product. There is no shortest connection from $(-1,0)$ to $(1,0)$ on $\R^{2,*}$.
A sequence of curves that is in $\R^{2,*}$ and converges to the curve $c:[-1,1]\rightarrow \R^2, t\mapsto (t,0)$
is readily constructed. Hence, the Riemannian distance between $(-1,0)$ and $(1,0)$ is $2$.
Yet, every curve connecting these points must go around the origin. The length-minimizing curve of length $2$ crosses the origin and is thus not an admissible curve on $\R^{2,*}$.}
Yet, candidates for shortest curves between points that are sufficiently close to each other can be obtained via a variational principle:
Given a parametric family of suitably regular curves $c_s:t\mapsto c_s(t)\in \mcM$, ${s\in(-\e,\e)}$ that connect
the same fixed endpoints $c_s(a) = p$ and  $c_s(b)=q$ for all $s$, one can consider the length functional
$s\mapsto L(c_s)$.
A curve $c = c_0$ is a first-order candidate for a shortest path between $p$ and $q$, if it is a critical point of the length functional, i.e., if $\frac{d}{ds}|_{s=0} L(c_s) = 0$. 
Such curves are called {\em geodesics}.
Differentiating the length functional leads to the so-called first variation formula \cite[\S 6]{Lee1997riemannian},
which, in turn, leads to the characterizing equation for geodesics:
\begin{definition}[Geodesics]
 \label{def:geodesic}
 A differentiable curve $c:[a,b]\rightarrow \mcM$ is called a {\em geodesic} (w.r.t. to a given Riemannian metric),
 if the {\em  covariant derivative} of its velocity vector field vanishes, i.e., 
 \begin{equation}
 \label{eq:geodesic}
  \frac{D\dot{c}}{dt}(t) = 0 \quad \forall t\in [a,b].
 \end{equation}
\end{definition}
\begin{remark}
 If a starting point $c(0) = p\in \mcM$ and a starting velocity $\dot{c}(0) = v\in T_p\mcM$ are specified, then
 the geodesic equation \eqref{eq:geodesic} translates to an initial value problem of second order
 with guaranteed existence and uniqueness of local solutions, \cite[p. 102]{AbsilMahonySepulchre2008}.
\end{remark}
An immediate consequence of \eqref{eq:geodesic} is that geodesics are constant-speed curves.
A formal introduction of the covariant derivative $\frac{D}{dt}$ along a curve is beyond the scope of this contribution,
and the interested reader is referred to, e.g., \cite[\S 4, \S 5]{Lee1997riemannian}.
To get some intuition, we introduce this concept for embedded Riemannian submanifolds $\mcM\subset \R^{n+d}$,
where the metric is the Euclidean metric of $\R^{n+d}$ restricted to the tangent bundle,
see also \cite[\S 20.12]{Gallier2011}:

A {\em vector field along a curve} $c:[a,b] \rightarrow \mcM$ is a differentiable map 
$v:[a,b]\rightarrow \R^{n+d}$ such that $v(t)\in T_{c(t)}\mcM$.
\footnote{The prime example for such a vector field is the curve's own velocity field $v(t) = \dot{c}(t)$.}
For every $p\in\mcM$, the ambient $\R^{n+d}$ decomposes into an orthogonal direct sum 
\[
 \R^{n+d} = T_p\mcM \oplus T_p\mcM^{\bot},
\]
where $T_p\mcM^{\bot}$ is the orthogonal complement of $T_p\mcM$ and orthogonality is w.r.t. the standard Euclidean inner product on $\R^{n+d}$.
Let $\Pi_p:\R^{n+d}\rightarrow T_p\mcM$ be the (base point-dependent) orthogonal projection onto the tangent space at $p$.
In this setting (and only in this), the covariant derivative of a vector field $v(t)$ along a curve $c(t)$ is the tangent component of $\dot{v}(t)$, i.e., 
$
 \frac{Dv}{dt}(t) = \Pi_{c(t)}(\dot{v}(t))
$.
As a consequence,
\begin{equation} 
\label{eq:GeodODE_euclid}
\frac{D\dot{c}}{dt}(t) = \Pi_{c(t)}(\ddot{c}(t)) 
\end{equation}
and the geodesics on Riemannian submanifolds with the metric induced by the ambient Euclidean inner product are precisely the constant-speed curves with acceleration vectors orthogonal to the corresponding tangent spaces,
i.e., $\ddot{c}(t) \in T_{c(t)}\mcM^{\bot}$.\\
{\bf Example:} On the unit sphere $S^2\subset \R^3$, the geodesics are great circles. When considered as curves in the ambient $\R^3$,
their acceleration vector points directly to the origin and is thus orthogonal to the corresponding tangent space.
When viewed as entities of $S^2$, these curves do not experience any acceleration at all.
\begin{center}
\begin{tikzpicture}[scale=.5]
\fill[black]  (3,0) circle (2pt) node [black,below left] {};
\draw[thick] (3,0) circle (2.5cm);
\draw[->]    (1.2322, 1.7678) -- (2.6464,    0.3536) node [midway, above]            {\hspace{0.9cm}$\ddot{c}(t)$};
\draw[->]    (1.2322, 1.7678) -- (-0.1820,    0.3536) node [midway, above]              {\hspace{-0.9cm}$\dot{c}(t)$};
\end{tikzpicture}
\end{center}
Mind that a constant-speed curve in $\R^n$ changes its direction only, when it experiences a non-zero acceleration.
In this sense, geodesics on manifolds are the counterparts to straight lines in the Euclidean space.

\tc{
In general, a covariant derivative, also known as a {\em linear connection},
is a bilinear mapping $(X,Y) \mapsto \nabla_XY$ that maps two vector fields $X,Y$ to a third vector field
$\nabla_XY$ in such a way that it can be interpreted as the directional derivative of $Y$ in the direction of $X$.
Of importance is the {\em Riemannian connection} or {\em Levi-Civita connection} that is compatible with a Riemannian metric \cite[Thm 5.3.1]{AbsilMahonySepulchre2008}, \cite[Thm 5.4]{Lee1997riemannian}. 
It is determined uniquely by the Koszul formula
\begin{eqnarray*}
   \label{eq:koszul}
   2\langle \nabla_XY, Z \rangle  &=& X(\langle Y, Z \rangle ) + Y(\langle Z,X \rangle) - Z(\langle X, Y \rangle)\\
   \nonumber
   && - \langle X, [Y,Z] \rangle - \langle Y, [X,Z] \rangle + \langle {Z}, {[X,Y]}\rangle
\end{eqnarray*}
and is used to define the {\em Riemannian curvature tensor}
$$
(X,Y,Z)\mapsto R(X,Y)Z=\nabla_X\nabla_YZ-\nabla_Y\nabla_XZ-\nabla_{[X,Y]}Z.\footnote{\tc{In these formulae, $[X,Y] = X(Y)-Y(X)$ is the Lie bracket of two vector fields.}}
$$
A Riemannian manifold is flat if and only if it is locally isometric to the Euclidean space, which holds if and only if the Riemannian curvature tensor vanishes identically \cite[Thm. 7.3]{Lee1997riemannian}. Hence, `flatness' depends on the Riemannian metric.
}

\subsection{Normal coordinates.}
\label{sec:normal_coords}
\begin{figure}[ht]
\begin{center}
\includegraphics[width=0.9\textwidth]{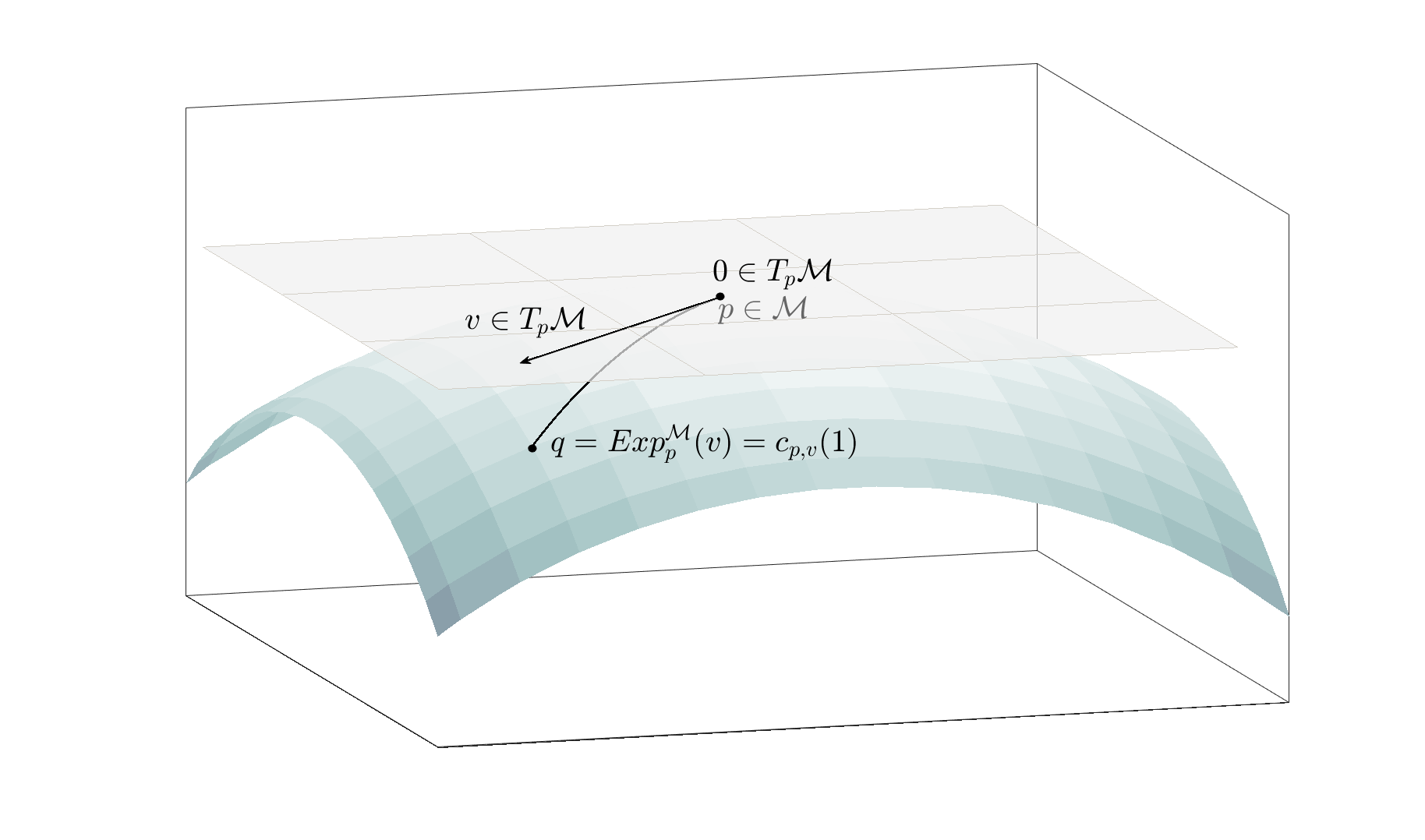}
\caption{The Riemannian exponential sends tangent vectors to end point of geodesic curves.
}
\label{fig:manifold_plotExp}
\end{center}
\end{figure}

The local uniqueness and existence of geodesics allows us to map a tangent vector $v\in T_p\mcM$
to the endpoint of a geodesic that starts from $p\in \mcM$ with velocity $v$.
Formalizing this principle gives rise to the {\em Riemannian exponential}
\begin{equation}
 \label{eq:Riemann_exp}
 \Exp^{\mcM}_p: T_p\mcM \supset B_{\e}(0) \rightarrow \mcM, \quad v\mapsto q:= \Exp^{\mcM}_p(v) := c_{p,v}(1).
\end{equation}
Here, $t\mapsto c_{p,v}(t)$ is the geodesic that starts from $p$ with velocity $v$ and  
$B_{\e}(0) \subset  T_p\mcM$ is the open ball with radius $\e$ and center $0$ in the tangent space\footnote{
For technical reasons, $\e>0$ must be chosen small enough such that $c_{p,v}(t)$ is defined on the unit interval $[0,1]$.},
see Fig. \ref{fig:manifold_plotExp}.
Note that we can restrict the considerations to unit-speed geodesics via 
$$\Exp^{\mcM}_p(v) := c_{p,v}(1) = c_{p, v/\|v\|}(t_v) = \Exp^{\mcM}_p\left(t_v \frac{v}{\|v\|}\right),$$ where $t_v = \|v\|$,
see \cite[\S 5., p. 72 ff.]{Lee1997riemannian} for the details.

For $\e>0$ small enough, the Riemannian exponential is a smooth diffeomorphism between $B_{\e}(0)$ and an open domain on $\mathcal{D}_p\subset\mcM$ around the point $p$. Hence, it is invertible. The smooth inverse map is called the {\em Riemannian logarithm} and is denoted by 
\begin{equation}
 \label{eq:Riemann_log}
\Log^{\mcM}_p:\mcM \supset \mathcal{D}_p \rightarrow B_\e(0)\subset T_p\mcM, \quad q \mapsto v:= (\Exp^{\mcM}_p)^{-1}(q),
\end{equation}
where $v$ satisfies $c_{p,v}(1) = q$.\\
Thus, the Riemannian logarithm is associated with the {\em geodesic endpoint problem}: Given $p,q\in \mcM$, find a geodesic that connects $p$ and $q$.
The Riemannian exponential map establishes a local parametrization of a small region around a location $p\in \mcM$
in terms of coordinates of the flat vector space $T_p\mcM$.
This is  referred to as representing the manifold in {\em normal coordinates}
\cite[\S III.8]{KobayashiNomizu1963}, \cite[Lem. 5.10]{Lee1997riemannian}.
Normal coordinates are radially isometric in the sense that the 
Riemannian distance between $p$ and $q = \Exp^\mcM_p(v)$ is exactly 
the same as the length of the tangent vector $\|v\|_p$ as
measured in the metric on $T_p\mcM$, provided that $v$ is contained in a
neighborhood of $0\in T_p\mcM$, where the exponential is invertible,
\cite[Lem. 5.10 \& Cor. 6.11]{Lee1997riemannian}.

Mind that \tc{the definition of the Riemannian exponential depends on the geodesics, which, in turn, depend on the chosen Riemannian metric -- via Definition \ref{def:RiemannMetric}.} Different metrics lead to different geodesics and thus to different exponential and logarithm maps.
%
%
%
%
%
\subsection{Matrix Lie groups and quotients by group actions}
\label{sec:LieGroups}
In general, a {\em Lie group} is a differentiable manifold $\mathcal{G}$ which
also has a group structure, such that the group operations
`multiplication' and `inversion',
$$\mathcal{G}\times
\mathcal{G}\ni(g, \tilde{g})\mapsto g\cdot \tilde{g}\in \mathcal{G}
\mbox{\hspace{0.1cm} and }\hspace{0.1cm} \mathcal{G}\ni g\mapsto g^{-1}\in \mathcal{G}
$$ 
are both smooth \cite{Gallier2011, Hall_Lie2015, godement2017introduction}.
A {\em matrix Lie group} $\mathcal{G}$ is a subgroup of $GL(n,\C)$ that is closed in $GL(n,\C)$.\footnote{but not necessarily in $\C^{n\times n}$.} This definition already implies that $\mathcal{G}$ is an embedded submanifold of $\C^{n\times n}$ \cite[Corollary 3.45]{Hall_Lie2015}.
Not all matrix groups are Lie groups and not all Lie groups are matrix Lie groups,
see \cite[\S 1.1 and \S 4.8]{Hall_Lie2015}. However, matrix Lie groups are arguably the most important class of Lie groups when it comes to practical applications and this exposition is restricted to this subclass.

Let $\mathcal{G}$ be an arbitrary matrix Lie group. \tc{When endowed with the bracket operator or {\em matrix commutator} $[V,W] = VW-WV$, the tangent space $T_I\mathcal{G}$} at the identity is called the {\em Lie algebra} associated with the Lie group $\mathcal{G}$, see \cite[\S 3]{Hall_Lie2015}. As such, it is denoted by $\mathfrak{g} = T_I\mathcal{G}$.
\tc{For any $A\in \mathcal{G}$, the function ``left-multiplication with $A$'' is a diffeomorphism
$L_A: \mathcal{G}\to \mathcal{G}, L_A(B) = AB$; its differential at a point $M\in \mathcal{G}$ is the isomporphism
$d(L_A)_M: T_M\mathcal{G}\to T_{L_A(M)}\mathcal{G}, d(L_A)_M(V) = AV$. Using this observation at $M=I$ shows that the}
tangent space at an arbitrary location $A\in \mathcal{G}$ is given by the translates (by left-multiplication) of the 
tangent space at the identity:
\begin{equation}
\label{eq:tangspaceshifts}
  T_A\mathcal{G} = T_{L_A(I)}\mathcal{G} =  A \mathfrak{g} = \left\{\Delta = AV\in \R^{n\times n}|\hspace{0.2cm} V\in \mathfrak{g}\right\},
\end{equation}
\cite[\S 5.6, p. 160]{godement2017introduction}.
The Lie algebra $\mathfrak{g} = T_I\mathcal{G}$ of $\mathcal{G}$ can equivalently be characterized as the set of all matrices $\Delta$ such that $\exp_m(t\Delta)\in \mathcal{G}$ for all $t\in \R$.
\tc{
The intuition behind this fact is that all tangent vectors are velocity vectors of smooth curves running on $\mathcal{G}$ (Definition \ref{def:tangentspace}) and that $c(t) = \exp_m(t\Delta)$ is a smooth curve starting from $c(0)=I$ with velocity $\dot{c}(0) =  \Delta$, see \cite[Def. 3.18 \& Cor. 3.46]{Hall_Lie2015} for the details.
}
\tc{By definition}, the exponential map\footnote{The exponential map of a Lie group must not be confused with the Riemannian exponential.} for a matrix Lie group is the matrix exponential restricted to the
corresponding Lie algebra, i.e. the tangent space at the identity $\mathfrak{g} = T_I\mathcal{G}$, \cite[\S 3.7]{Hall_Lie2015},
\[
 \exp_m|_{\mathfrak{g}}: \mathfrak{g}\rightarrow \mathcal{G}.
\]
In general, a Lie algebra is a vector space with a linear, skew-symmetric bracket operation, called {\em Lie bracket}
$[\cdot,\cdot]$ that satisfies the {\em Jacobi identity}.
\[
 [X, [Y,Z]] + [Z,[X,Y]] + [Y,[Z,X]] = 0.
\]
\paragraph*{\tc{Quotients of Lie groups by closed subgroups}}
\tc{
In many settings, it is important or sometimes even necessary to consider certain points $p,q$ on a given differentiable manifold $\mcM$ as equivalent.
Consider the following example.
 \begin{example}
 \label{ex:QuotBasic}
  Let $U\in \R^{n\times \p}$ feature orthonormal columns so that $U^TU=I_\p$. We may extend the columns of $U=(u^1,\ldots,u^\p)$ to an orthogonal matrix $Q= (u^1,\ldots,u^\p, u^{\p+1}, \ldots,u^n)\in O(n)$.
  Let $I_\p\times O(n-\p):= \left\{\begin{pmatrix}
                            I_\p &0\\
                             0    & R
                           \end{pmatrix}| \hspace{0.1cm} R\in O(n-\p)\right\}$.                           
  This is actually a closed subgroup of $O(n)$, 
  in symbols $(I_\p\times O(n-\p))\leq O(n)$.
  The action $\tilde{Q} = Q\Phi$ with any orthogonal matrix $\Phi\in I_\p\times O(n-\p)$ 
  preserves the first $\p$ columns of $Q$. Hence, we may identify $U$ with the equivalence class 
  $[Q] = \{Q\Phi| \Phi \in I_\p\times O(n-\p)\}\subset O(n)$.
  In Sections \ref{sec:Stiefel} and \ref{sec:Grassmann}, we will see that this example establishes the Stiefel manifold of ONBs and eventually also the Grassmann manifold of subspaces as quotients of the orthogonal group $O(n)$.
 \end{example}
Note that in the example, the equivalence relation is induced by actions of the Lie group $I_\p\times O(n-\p)$. 
Quotients that arise from such group actions are important examples of {\em quotient manifolds}.
The following Theorems \ref{thm:QuotMan} and \ref{thm:HomSpaceConst} cover this example as well as all other cases
of quotient manifolds that are featured in this work.
First, group actions need to be formalized.}
\tc{
\begin{definition}
  \label{def:orbitspace}(cf. \cite[p. 162,163]{Lee2012smooth})
  Let $\mathcal{G}$ be a Lie group, $\mcM$ be a smooth manifold, and let $\mathcal{G}\times \mcM \to \mcM, (g,p)\mapsto g\cdot p$ be a
  {\em left action of $\mathcal{G}$ on $\mcM$.}\footnote{\tc{The theory for right actions is analogous.
  In all cases considered in this work, $\mcM$ is a matrix manifold so that ``$\cdot$'' is the usual matrix product.  }}
  The {\em orbit relation} on $\mcM$ induced by $\mathcal{G}$ is defined by
   \[
    p\simeq q :\Leftrightarrow \exists g\in \mathcal{G}: \quad g\cdot p = q.
   \]
   The equivalence classes are the {\em $\mathcal{G}$-orbits} $[p]:=\mathcal{G}p := \{g\cdot p|\hspace{0.1cm} g\in \mathcal{G}\}.$
 The {\em orbit space} is denoted by $\mcM/\mathcal{G}:= \{[p]|\hspace{0.1cm} p\in \mcM\}$.
 The {\em quotient map} sends a point to its $\mathcal{G}$-orbit via
 $\Pi : \mcM\to \mcM/\mathcal{G}, \hspace{0.1cm} p\mapsto [p]$.
 The action is {\em free}, if every isotropy group
$\mathcal{G}_p:= \{g\in \mathcal{G}|\hspace{0.1cm} g\cdot p = p\}$ is trivial, $\mathcal{G}_p=\{e\}$.
 \end{definition}
\begin{theorem}
\label{thm:QuotMan}
(Quotient Manifold Theorem, cf. \cite[Thm. 21.10]{Lee2012smooth})
Suppose $\mathcal{G}$ is a Lie group acting
{\em smoothly, freely, and properly} on a smooth manifold $\mcM$. 
Then the orbit space $\mcM/\mathcal{G}$ is a manifold of dimension $\dim\mcM - \dim \mathcal{G}$,
and has a unique smooth structure such that the quotient map 
$\Pi : \mcM\to \mcM/\mathcal{G}, p\mapsto [p]$
is a smooth submersion.\footnote{\tc{i.e. a smooth surjective mapping such that the differential is surjective at every point.}}
In this context, $\mcM$ is called the total space and $\mcM/\mathcal{G}$ is the quotient (space).
\end{theorem}
}
\tc{
A special case is Lie groups under actions of Lie subgroups.
 \begin{definition}{\cite[\S 21, p. 551]{Lee2012smooth}}
  \label{def:coset}
  Let $\mathcal{G}$ be a Lie group and $\mathcal{H}\leq \mathcal{G}$ be a Lie subgroup. 
 For $g\in \mathcal{G}$, a subset of $\mathcal{G}$ of the form 
   $
    [g]:= g\mathcal{H} = \{g\cdot h|\hspace{0.1cm} h\in \mathcal{H}\}
   $
  is called a {\em left coset of $\mathcal{H}$}.
  The left cosets form a partition of $\mathcal{G}$,
  and the quotient space determined by this partition is called the {\em left coset space of $\mathcal{G}$ modulo $\mathcal{H}$}, 
  and is denoted by $\mathcal{G}/\mathcal{H}$.
 \end{definition}
Coset spaces of Lie groups are again smooth manifolds:
 \begin{theorem}
  \label{thm:HomSpaceConst}
  (cf. \cite[Thm 21.17, p. 551]{Lee2012smooth})
  Let $\mathcal{G}$ be a Lie group and let $\mathcal{H}$ be a closed subgroup of $\mathcal{G}$. The left coset space $\mathcal{G}/\mathcal{H}$ is a manifold of dimension $\dim \mathcal{G} - \dim \mathcal{H}$ with a unique differentiable structure such that the quotient map 
  $\Pi: \mathcal{G} \to  \mathcal{G}/\mathcal{H}, g\mapsto [g]$ is a smooth submersion.
 \end{theorem}
}
\tc{
In general, if $\pi: \mcM\to \mcN$ is a surjective submersion between two manifolds $\mcM$ and $\mcN$,
then for any $q\in \mcN$, the the preimage $\pi^{-1}(q) \subset \mcM$ is called the {\em fiber over $q$}, and is denoted  by $\mcM_q$. Each fiber $\mcM_q$ is itself a closed, embedded submanifold by the inverse function theorem.
If $\mcM$ has a Riemannian metric $\langle \cdot,\cdot\rangle_p^\mcM$, 
then at each point $p\in \mcM$, the tangent space $T_p\mcM$
decomposes into an {\em orthogonal  direct sum} $T_p\mcM = T_p\mcM_{\pi(p)}\oplus (T_p\mcM_{\pi(p)})^\bot$.
The tangent space of the fiber $T_p\mcM_{\pi(p)}=: V_p$ is the called the {\em vertical space},
its orthogonal complement $H_p := V_p^{\bot}$ is the {\em horizontal space}.
The vertical space is the kernel $V_p = \ker(d\pi_p)$ of the differential 
$d\pi_p: T_p\mcM \to T_{\pi(p)}\mcN$;
the horizontal space is isomorphic to $T_{\pi(p)}\mcN$.
This allows to identify $H_p \cong T_{\pi(p)}\mcN$, see \cite[Fig. 3.8., p. 44]{AbsilMahonySepulchre2008}
for an illustration. This construction helps to compute tangent spaces of quotients, if the tangent space of the total space is known.
}


 If $\mcN$ is a quotient as in Theorem \ref{thm:QuotMan} or Theorem \ref{thm:HomSpaceConst} and if $\Pi:\mcM\to \mathcal{N}$ is the corresponding quotient map, then $\Pi$ can be turned into a {\em Riemannian submersion}, i.e., a submersion that is compatible with the Riemannian metric in the sense that $d\Pi$ preserves inner products of horizontal vectors, see
 \cite[Chap. 8, Sec. 5, ex. 8.-9.]{DoCarmo2013riemannian}.
For every tangent vector $w\in T_{\Pi(p)}\mcN$ there is $\bar x =\bar v+\bar w\in V_p \oplus H_p = T_p\mcM$ such that 
$d\Pi_p(\bar x) = w$. The horizontal component $\bar w$ is unique and is called the {\em horizontal lift} of $w$.
By relying on horizontal lifts, a Riemannian metric on the quotient can be defined by
\begin{equation}
	\label{eq:quotMetric}
	\langle w_1,w_2 \rangle^{\mcN}_{\Pi(p)} 
	:= \langle \bar w_1, \bar w_2 \rangle^{\mcM}_p 
\end{equation}
for $w_1,w_2 \in T_{\Pi(p)}\mcN$. With respect to this (and only this) metric, the quotient map is a {\em local isometry} between the horizontal space $H_p$ and  $T_{\Pi(p)}\mcN$.
As a consequence, horizontal geodesics in $\mcM$ are mapped to geodesics in $\mcN$ under $\Pi$. 
Horizontal geodesics are geodesics in the total space, whose velocity field stays in the horizontal space for all time $t$.

Theorem \ref{thm:HomSpaceConst} additionally establishes $\mathcal{G}/\mathcal{H}$ as a
{\em homogeneous  space}, i.e. a smooth manifold $\mcM$ endowed with a 
  {\em transitive smooth action} by a Lie group (cf. \cite[\S 21, p. 550]{Lee2012smooth}).
In the setting of the theorem, the group action is given by
the left action of $\mathcal{G}$ on $\mathcal{G}/\mathcal{H}$ given by $g_1\cdot [g_2] := [g_1\cdot g_2]$.
A transitive action allows us to transport a location $p\in \mcM$ to any other location $q\in \mcM$.

%
%
\section{Interpolation on non-flat manifolds}
\label{sec:ManInterp}
%
%
%
When working with matrix manifolds, the data is usually given in
extrinsic coordinates, see Section \ref{sec:diffgeo_intro}. For example, data on
the compact Stiefel manifold $St(n,\p)= \{U\in\R^{n\times \p}| U^TU=I_\p\}$, $\p\leq n$, 
is given in form of $n$-by-$\p$ matrices.
These matrices feature $n\p$ entries while the intrinsic number of degrees of freedom, i.e., the intrinsic dimension is turns out to be $n\p-\frac{1}{2}\p(\p+1)$, see Section \ref{sec:Stiefel}.
Essentially, the practical obstacle associated with data interpolation on matrix manifolds arises from this fact. Given, say, $k$ matrices on $St(n,\p)$ in extrinsic coordinates, interpolating entry-by-entry will most certainly lead to interpolants that do not feature orthogonal columns and thus are {\em not} points on the Stiefel manifold.
Likewise, entry-by-entry interpolation of positive definite matrices is not guaranteed to produce another positive definite matrix.

There are essentially two different approaches to address this issue:
Performing the interpolation on the tangent space of the manifold and 
using the Riemannian barycenter or Riemannian center of mass as an interpolant.
Both will be explained in more detail in the next two subsections.\footnote{German speaking readers may find an introduction that addresses a general scientific audience in \cite{SanderRundbrief2015}.}
\subsection{Interpolation in normal coordinates}
\label{sec:tangInt}
As outlined in Section \ref{sec:diffgeo_intro}, every location $p\in\mathcal{M}$
on an $n$-dimensional differentiable manifold features a small neighborhood $\mathcal{D}_p$
that is the domain 
of a coordinate chart $x:\mathcal{M}\supset\mathcal{D}_p\rightarrow D_{x(p)}  \subset \R^n$
that maps bijectively onto an open set $D_{x(p)}\subset \R^n$.
Therefore, for a sample data set $\{p_1,\ldots, p_k\}\subset \mathcal{D}_p$ that is completely
contained in the domain of a single coordinate chart $x$, interpolation can be performed as follows:
\begin{enumerate}
 \item Map the data set to $D_{x(p)}$: Calculate $v_1 = x(p_1),\ldots, v_k = x(p_k)\in D_{x(p)}$.
 \item Interpolate in $D_{x(p)}$ to produce the interpolant $v^* \in D_{x(p)}$.
 \item Map back to manifold: compute $p^* = x^{-1}(v^*)\in \mathcal{D}_p$.
\end{enumerate}
In principle, any coordinate chart may be applied.
In practice, the challenge is to find a suitable coordinate chart that can be evaluated efficiently.
Moreover, it is desirable that the chosen chart preserves the geometry of the original data set as well as possible.\footnote{There are no isometric coordinate charts on a non-flat manifold, see \cite[Thm 7.3]{Lee1997riemannian}.}
The standard choice is to use {\em normal coordinates} as introduced in Section \ref{sec:normal_coords}.
This means that the Riemannian logarithm is used as the coordinate chart
%
\[
  \Log^{\mcM}_p:\mcM \supset \mathcal{D}_p \rightarrow B_\e(0)\subset T_p\mcM
\]
with the Riemannian exponential 
\[
 \Exp^{\mcM}_p: T_p\mcM \supset B_{\e}(0) \rightarrow \mathcal{D}_p\subset\mcM
\]
as the corresponding parameterization.
The general procedure of data interpolation via the tangent space is formulated as 
Algorithm \ref{alg:TangInt}.
\begin{algorithm}
\caption{Interpolation in normal coordinates.}
\label{alg:TangInt}
\begin{algorithmic}[1]
  \REQUIRE{Data set $\{p_1,\ldots, p_k\}\subset \mcM$.}
  \STATE{Choose $p_{i}\in \{p_1,\ldots, p_k\}$ as a base point.}
  \STATE{Check that $\Log^{\mcM}_{p_{i}}(p_j)$ is well-defined for all $j=1,\ldots,k$.}
  \FOR{$j=1,\ldots,k$}
    \STATE{Compute $v_j := \Log^{\mcM}_{p_{i}}(p_j) \in T_{p}\mcM$.}
  \ENDFOR
  \STATE{Compute $v^*$ via Euclidean interpolation of $\{v_1,\ldots, v_k\}$.}
  \STATE{Compute $p^* := \Exp^{\mcM}_{p_{i}}(v^*)$}
  \ENSURE{$p^*\in \mcM$.}
\end{algorithmic}
\end{algorithm}
\begin{remark}
\label{rem:onInterpNormCoord}
There are a few facts that the practitioner needs to be aware of:
 \begin{enumerate}
  \item The interpolation procedure of Algorithm \ref{alg:TangInt} depends on which sample point is selected to act as the base point. Different choices may lead to different interpolants.\footnote{In the practical applications considered in \cite{AmsallemFarhat2008}, it was observed that the base point selection has only a minor impact on the final result.}
  \item For {\em matrix manifolds}, the tangent space is often also given in {\em extrinsic coordinates}. This means that an entry-by-entry interpolation of the matrices that represent the tangent vectors may lead to an interpolant that is not in the tangent space. 
  \tc{As an illustrative example, consider the Grassmannian $Gr(n,\p)$.
  Matrices} $\Delta_1,\ldots,\Delta_k\in T_{[U]}Gr(n,\p)$ are characterized by $U^T\Delta_j = 0$.
  Entry-by-entry interpolation in the tangent space may potentially result in a matrix $\Delta^*$ that is not orthogonal to the base point $U$, i.e. $U^T\Delta^*\neq 0$,  see \cite[\S 2.4]{Zimmermann2014}.
  
  \tc{In general, because of the vector space structure of the tangent space of any manifold $\mcM$, it is sufficient to use an interpolation method that expresses the interpolant in $T_p\mcM$
  as a weighted linear combination of the sampled tangent vectors
  $v_1,\ldots,v_k\in T_{p}\mcM$
  \[
    v^* = \sum_{j=1}^k \omega_j v_j.
  \]
  }
  Amongst others, {\em linear interpolation}, {\em Lagrange and Hermite interpolation}, {\em spline interpolation} and interpolation via {\em radial basis functions} fulfill this requirement. As an aside, the interpolation procedure is computationally less expensive, since it works on the weight coefficients $\omega_j$ rather than on every single entry.
 \end{enumerate}
 \label{rem:tangInterp}
\end{remark}
\paragraph*{Quasi-linear interpolation of trajectories via geodesics}
In this paragraph, we address applications, where the sampled manifold data features a univariate parametric dependency.
The setting is as follows.
Let $\mcM$ be a Riemannian manifold and suppose that there is a trajectory
\[
 c:[a,b] \rightarrow \mcM,\quad \mu \mapsto c(\mu) 
\]
on $\mcM$ that is sampled at $k$ instants $\mu_1,\ldots,\mu_k\in [a,b]$.
Then, an interpolant $\hat{c}$ for $c$ can be computed via Algorithm \ref{alg:GeoInt}.
\begin{algorithm}
\caption{Geodesic interpolation}
\label{alg:GeoInt}
\begin{algorithmic}[1]
  \REQUIRE{Data set $\{c(\mu_1),\ldots, c(\mu_k)\}\subset \mcM$ sampled from a curve $c:\mu\rightarrow c(\mu)$,
  unsampled instant $\mu^*\in [\mu_j, \mu_{j+1}]$.}
 
  \STATE{Compute $v_{j+1} := \Log^{\mcM}_{c(\mu_{j})}(c(\mu_{j+1})) \in T_{c(\mu_{j})}\mcM$.}
  \STATE{Compute $\hat{c}(\mu^*) := \Exp^{\mcM}_{c(\mu_{j})}\left(\frac{\mu^*-\mu_j}{\mu_{j+1}-\mu_j} v_{j+1}\right)$}
  \ENSURE{$\hat{c}(\mu^*)\in \mcM$ interpolant of $c(\mu^*)$.}
\end{algorithmic}
\end{algorithm}
The interpolants at $\mu\in [\mu_j,\mu_{j+1}]$ that are  output by Algorithm \ref{alg:GeoInt} lie on the unique geodesic connection between
the points $c(\mu_j)$ and $c(\mu_{j+1})$. Hence, it is the straightforward manifold analogue of linear interpolation \tc{ and is base-point independent}.

\begin{figure}[ht]
\centering
\includegraphics[width=1.0\textwidth]{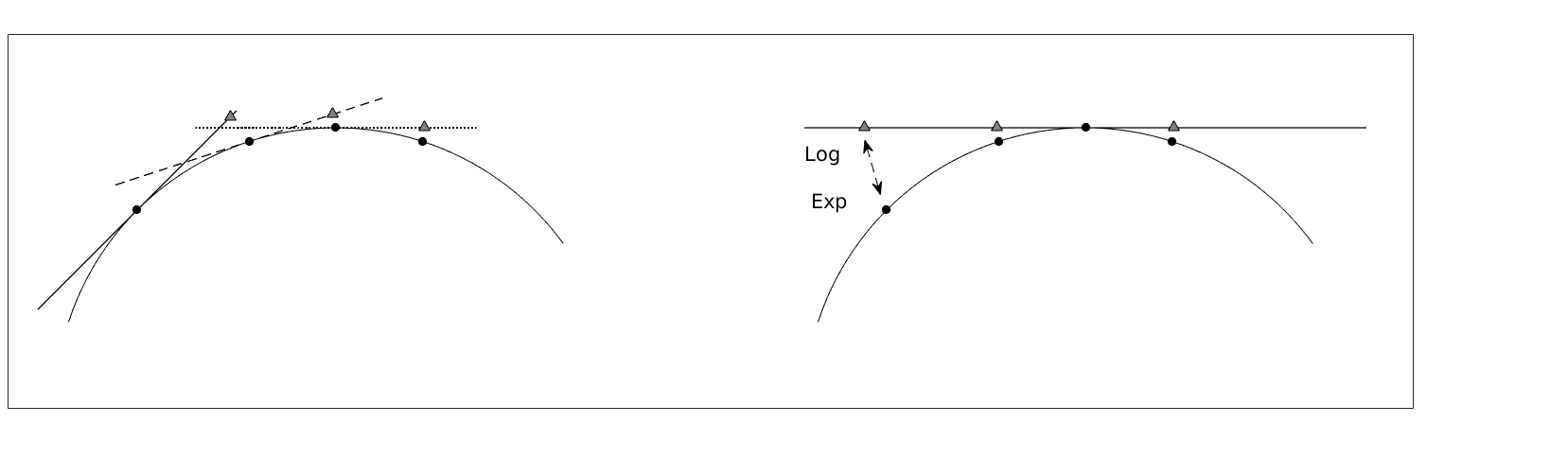}
\caption{\tc{Illustration of the course of action of Algorithms \ref{alg:TangInt} and \ref{alg:GeoInt}. Algorithm \ref{alg:TangInt} (right) first maps all data points to a selected fixed tangent space. In Algorithm \ref{alg:GeoInt} (left), two points $p_j=c(\mu_j)$ and $p_{j+1}=c(\mu_{j+1})$ are connected by a geodesic line, then the base is shifted to point $p_{j+1}$ and the procedure is repeated. 
}}
\label{fig:plot_InterpPrinciple}
\end{figure}

\tc{
The generic formulation of Algorithm \ref{alg:TangInt} allows to employ higher-order interpolation methods. However, this does not necessarily lead to more accurate results: the overall error depends not only on the interpolation error within the tangent space but also on the distortion caused by mapping the data to a selected (fixed) tangent space, see Fig. \ref{fig:plot_InterpPrinciple}.
}
%
%
%

Algorithms \ref{alg:TangInt} and \ref{alg:GeoInt} can be applied in practical applications,
where the Riemannian exponential and logarithm mappings are known in explicit form.
Applications in parametric model reduction that consider matrix manifolds include \cite{Degroote_etal2010} ($GL(n)$-data),
\cite{AmsallemFarhat2008, Nguyen2012, Zimmermann2014} (Grassmann-data), \cite{Zimmermann2018} (Stiefel data)
and \cite{AmsallemFarhat2011, Rahman_etal2005} ($SPD(n)$-data).
\subsection{Interpolation via the Riemannian center of mass}
\label{sec:CenterOfMass:Int}
As pointed out in Remark \ref{rem:tangInterp}, interpolation of manifold data via the back and forth mapping of a complete data set of sample points between the manifold and its tangent space depends on the chosen base point. As a consequence, sample points may experience an uneven distortion under the projection onto the tangent space, see Fig. \ref{fig:plot_InterpPrinciple} (right).
An approach that avoids this issue is to interpret interpolation as the task of finding suitably weighted Riemannian centers of mass.
This concept was introduced in the context of geodesic finite elements in \cite{Sander:2016:GFE, Grohs_quasiRiemann2013}.

The idea is as follows:
The Riemannian center of mass\footnote{Here, we introduce this for discrete data sets; for centers w.r.t. a general mass distribution, see the original paper \cite{Karcher1977}, Section 1.} or Fr{\'e}chet mean of a sample data set
$\{p_1,\ldots,p_k\}\in \mcM$ on a manifold
with respect to the scalar weights $w_i\geq 0$, $\sum_{i=0}^k w_i = 1$
is defined as the minimizer(s) of the Riemannian objective function
\[
 \mcM\ni q \mapsto f(q) = \frac{1}{2} \sum_{i=1}^k w_i \dist(q, p_i)^2,
\]
where $\dist(q,p_i)$ is the Riemannian distance of \eqref{eq:RiemannDist}.
This definition generalizes the notion of the barycentric mean in Euclidean spaces.
However, on curved manifolds, the global center might not be unique. Moreover,
local minimizers may appear. For more details, see \cite{Karcher1977} and \cite{AfsariTronVidal2013}, which also give uniqueness criteria.\\
Interpolation is now performed by computing weighted Riemannian centers.
More precisely, 
let $\mu_1,\ldots,\mu_k\subset \R^d$ be sampled parameter locations and let 
$p_i=p(\mu_i)\in \mcM$, $i=1,\ldots,k$ be the corresponding sample locations on $\mcM$.
Interpolation is within the convex hull $\text{conv}\{\mu_1,\ldots,\mu_k\}\subset \R^d$ of the samples.

Let $\{\varphi_i:\mu\mapsto \varphi_i(\mu)| i=1,\ldots,k\}$ be a suitable set of interpolation functions with $\varphi_i(\mu_j) = \delta_{ij}$, $\sum_i \varphi_i(\mu)\equiv 1$, say Lagrangians \cite{Sander:2016:GFE}, splines \cite{Grohs_quasiRiemann2013} or radial basis functions \cite{Buhmann2003}.
Then, the interpolant $p^* \approx p(\mu^*)\in \mcM$ at an unsampled parameter 
location $\mu^*\in \text{conv}\{\mu_1,\ldots,\mu_k\}$ is defined 
as the minimizer of
\begin{equation}
 \label{eq:baryInterp}
 p^* = \argmin_{q\in \mcM} f(q) = \frac12 \sum_{i=1}^k \varphi_i(\mu^*) \dist(q,p_i)^2.
\end{equation}
At a sample location $\mu_j$, one has indeed that
$$\sum_{i=1}^k \varphi_i(\mu_j) \dist(q,p_i)^2
= \sum_{i=1}^k \delta_{ij} \dist(q,p_i)^2= \dist(q,p_j)^2,
$$
which has the
unique global minimum at $q = p_j$.

Computing $p^*$ requires to solve a Riemannian optimization problem.
The simplest approach is a gradient descent method \cite{AfsariTronVidal2013, AbsilMahonySepulchre2008}.
The gradient of the objective function $f$ in \eqref{eq:baryInterp} is
\begin{equation}
 \label{eq:baryGrad}
 \nabla f_q = - \sum_{i=1}^k \varphi_i(\mu^*) \Log^\mcM_q(p_i) \in T_q\mcM.
\end{equation}
see \cite[Thm 1.2]{Karcher1977}, \cite[\S 2.1.5]{AfsariTronVidal2013},
\cite[eq. (2.4)]{Sander:2016:GFE}. Hence, just like interpolation in the tangent space, the interpolation via the Riemannian center can be pursued only in applications, where the Riemannian logarithm can be computed.
A generic gradient descent algorithm to compute the barycentric interpolant
for a function $p:\R^d \ni \mu \mapsto p(\mu) \in \mcM$ reads as follows.
\begin{algorithm}
\caption{Interpolation via the weighted Riemannian center \cite{Rentmeesters2011,AfsariTronVidal2013}.}
\label{alg:BaryInt}
\begin{algorithmic}[1]
  \REQUIRE{Sample data set $\{p_1=p(\mu_1),\ldots, p_k=p(\mu_k)\}\subset \mcM$,
  unsampled parameter location $\mu^*\in \text{conv}(\mu_1,\ldots, \mu_{k})\subset \R^d$,
  initial guess $q_0$, convergence threshold $\tau$.}
  \STATE{$k:=0$}
  \STATE{Compute $\nabla f_{q_k }$ according to \eqref{eq:baryGrad}}
  \WHILE{$\|\nabla f_{q_k }\|_{q}> \tau$}
  \STATE{select a step size $\alpha_k$}
  \STATE{$q_{k+1} := \Exp^{\mcM}_{q_k}\left(-\alpha_k \nabla f_{q_k }\right)$}
  \STATE{$k:= k+1$}
  \ENDWHILE
  \ENSURE{$p^*:=q_k\in \mcM$ interpolant of $p(\mu^*)$.}
\end{algorithmic}
\end{algorithm}
An implementation of this (type of) method for finding the Karcher mean in $SO(3)$ is discussed in \cite{Rentmeesters2011}. Of course, Riemannian analogues to more sophisticated nonlinear optimization methods may also be employed, see \cite{AbsilMahonySepulchre2008}.

In the context of model reduction, the benefits of interpolation via weighted Riemannian centers and the computational costs of solving the associated Riemannian optimization problem  must be juxtaposed.
\subsection{Additional approaches}
A large variety of sophistications and further manifold interpolation techniques exists in the literature:
The acceleration-minimi\-zing property of cubic splines in the Euclidean space 
can be generalized to Riemannian manifolds in form of a variational problem
\cite{Noakes1989, Crouch1995, CAMARINHA2001107, Steinke2008,BOUMAL2011,Samir2012, Kim2018},
see also \cite{Noakes2007} and references therein.
%
Moreover, the construction concepts of B{\'e}zier curves and the De Casteljau-algorithm \cite{bartels1995introduction}
can be transferred to Riemannian manifolds \cite{Noakes2007, Krakowski2015SOLVINGIP, Polthier2013, AbsGouseWirth2016, SAMIR2019}.
B{\'e}zier curves in Euclidean spaces are polynomial splines that rely on a number of so-called control points.
To obtain the value of a B{\'e}zier curve at time $t$, a recursive sequence of straight-line convex combinations
between pairs of control points must be computed.
The transition of this technique to Riemannian manifolds is via replacing the inherent straight lines with geodesics \cite{Noakes2007}.
Another option is to conduct the B{\'e}zier/De Casteljau-algorithm in the tangent space 
and to transfer the results to the manifold via a geodesic averaging of the spline arcs 
that were constructed in the tangent spaces at the first and the last control point, respectively, see \cite{GouseMassartAbsil2018}.

Derivative information may also be  incorporated in inter\-po\-lation schemes on Riemannian manifolds.
A Hermite-type method that is specifically tailored for interpolation problems on the Grassmann manifold is sketched in \cite[\S 3.7.4]{Amsallem2010}.
General Hermitian manifold interpolation in compact, connected Lie groups with a bi-invariant metric has been considered in \cite{Jakubiak2006}.
A practical approach to conduct first-order Hermite interpolation of data on arbitrary Riemannian manifolds is discussed in \cite{Zimmermann_2019}.
%
\subsection{Quasi-linear extrapolation on matrix manifolds}
\label{sec:extrapol}
%
%
In application scenarios, where both snapshot data of the full-order model and derivative information are at hand, various approaches have been suggested 
to exploit the latter.
On the one hand, derivatives can be used for improving the ROMs accuracy and approximation quality by constructing POD bases that incorporate snapshots and snapshot derivatives \cite{CPOD_Carlberg_2011,HinzeVolkwein2005, Ito, Zimmermann2012}.
On the other hand, snapshot derivatives enable
to parameterize the ROM bases and subspaces or to perform sensitivity analyses \cite{wem, HayBorggaardPelletier2009,HayBorggaardAkhtarPelletier2010, Zimmermann2016}.
In this section, we outline an approach to transfer the idea of extrapolation and parameterization via local linearizations to manifold-valued functions.
The underlying idea is comparable to the trajectory piece-wise linear (TPWL) method \cite{Rewienski2004}. 
Yet, TPWL linearizes the full-order model prior to the ROM projection, whereas here, we consider linearizing ROM building blocks like the reduced orthogonal bases, reduced subspaces or reduced system matrices.
\paragraph*{A geometric first-order Taylor approximation}
Any differentiable function $f:\R^n\to\R^n$ can be linearized via a first-order Taylor expansion.
A step ahead of size $t$ in direction $d\in \R^n$ gives
$f(x_0+td) = f(x_0) + tDf_{x_0}(d) + \mathcal{O}(t^2).$
When considering $t\mapsto c(t) := f(x_0+td)$ as a curve, then the 
first-order Taylor approximant is the straight line $g:t\mapsto c(0) + \dot c(0) t$.
Such first order linearization often serves for extrapolating a given nonlinear function in a neighborhood of a selected expansion point. For doing so, the starting point $c(0)$ and the starting velocity $\dot c(0)$ must be available.
This procedure translates to the manifold setting, when straight lines are replaced with geodesics.

Let $\mu\in\R$ be a scalar parameter and let $c:\mu\mapsto c(\mu)\in \mcM$ be a curve on a submanifold $\mcM$.
For given initial values $c(\mu_0) = p_0\in \mcM$ and $\dot c (\mu_0) = v_0\in T_{p_0}\mcM$,
the corresponding unique geodesic $c_{p_0,v_0}$ is expressed via the Riemannian exponential as
$$
   c_{p_0,v_0}: \mu \rightarrow \mathcal{M}, \quad \mu\mapsto  \Exp_{p_0}^\mcM(\mu v_0).
$$
%
%
%
%
%
%
%
\begin{algorithm}
\caption{Geodesic extrapolation.}
\label{alg:GeoExt}
\begin{algorithmic}[1]
  \REQUIRE{Scalar parameter $\mu_0\in\R$, initial values $c(\mu_0)\in \mcM, \dot c(\mu_0)\in T_{c(\mu_0)}\mcM$ sampled from a curve $c:\mu\rightarrow c(\mu)\in \mcM$, parameter value $\mu^* >0$.}
  \STATE{Compute $\hat{c}(\mu_0 + \mu^*) := \Exp^{\mcM}_{c(\mu_0)}\left(\mu^* \dot c(\mu_0)\right)$}
  \ENSURE{$\hat{c}(\mu_0 + \mu^*)\in \mcM$ extrapolant of $c(\mu_0 + \mu^*)$.}
\end{algorithmic}
\end{algorithm}
%
%
%
\paragraph*{Example: Extrapolating POD basis matrices}
As outlined in Section \ref{sec:introInterp4pMOR},
snapshot POD works by collecting state vector snapshots,
$x^1:=x(t_1,\mu_0),...,x^{m}:=x(t_m,\mu_0)\} \in \R^n$
followed by an SVD of the snapshot matrix $\left(x^1,...,x^{m}\right)(\mu_0)=: \Sn(\mu_0) = \U(\mu_0)\Sigma(\mu_0) \Z^T(\mu_0)$.
Here, the matrix dimensions are $\U(\mu_0)\in \R^{n\times m}$, $\Sigma(\mu_0)\in\R^{m\times m}$, 
$\Z(\mu_0)\in \R^{m\times m}$.
The objective is to approximate $\U(\mu_0 + \mu)$ for a small $\mu>0$ based on the data $\U(\mu_0), \dot \U(\mu_0)$,
\tc{where $\U(\mu_0)$ is a point on the Stiefel manifold $St(n,m)$ and $\dot \U(\mu_0)$ is a tangent vector, see Section \ref{sec:Stiefel_intro_and_numerics}}.
{\em Differentiating the SVD.}
If the snapshot matrix function $\mu \mapsto \Sn(\mu)\in \R^{n\times m}$ 
is smooth in the neighborhood of $\mu_0\in \R$ and 
if the singular values of $\Sn(\mu_0)$ are mutually distinct\footnote{This condition can be relaxed, see the results of \cite[\S 7]{Alekseevsky1998}.}, then
the singular values and both the left and the right singular vectors 
are differentiable in $\mu \in [\mu_0 -\delta\mu, \mu_0+\delta\mu]$
for $\delta\mu$ small enough.
For brevity, let $\dot{\Sn} = \frac{d\Sn}{d\mu}(\mu_0)$ denote the derivative
with respect to $\mu$ evaluated in $\mu_0$ and so forth.
Let $\mu \mapsto \Sn(\mu) = \U(\mu)\Sigma(\mu) \Z(\mu)^T \in \R^{n\times m}$ and let $C(\mu) = (\Sn^T\Sn)(\mu)$. 
Let $u^j$ and $v^j$, $j=1,\ldots,m$ denote the columns of $\U(\mu_0)$ and $\Z(\mu_0)$, respectively.
It holds
\begin{eqnarray}
 \label{eq:SVD-diff1}
        \dot{\sigma}_j &=& (u^j)^T \dot{\Sn} v^j, (j = 1,\ldots, m),\\
 \label{eq:SVD-diff2}
        \dot{\Z} &=& \Z A, \mbox{ where } A_{ij} =
            \left\{
                \begin{array}{ll}
                 \frac{\sigma_j(u^j)^T\dot\Sn v^i + \sigma_i(u^i)^T\dot\Sn v^j}{(\sigma_j + \sigma_i)(\sigma_j - \sigma_i)}, & i\neq j\\
                0,          & i=j
                \end{array}
            \right. (i,j = 1,\ldots, m),\\
 \label{eq:SVD-diff3}
    \dot{\U} &=& \dot{\Sn}\Z\Sigma^{-1} + \Sn\dot{\Z}\Sigma^{-1} + \Sn\Z \dot{\Sigma}^{-1}
= \left(\dot{\Sn} \Z + \U(\Sigma A - \dot{\Sigma}) \right) \Sigma^{-1}.
\end{eqnarray}
%
A proof can be found in \cite{HayBorggaardPelletier2009}. 
Note that $\U^T(\mu_0)\dot{\U}(\mu_0)$ is skew-symmetric so that indeed $\dot{\U}(\mu_0)=:\Delta(\mu_0) \in T_{\U(\mu_0)}St(n,m)$.
The above equations hold in approximative form for the truncated SVD.
For convenience, assume that $\U(\mu_0)\in St(n,\p)$ is now the truncated to $\p\leq m$ columns.

\begin{figure}[ht]
\centering
\includegraphics[width=0.7\textwidth]{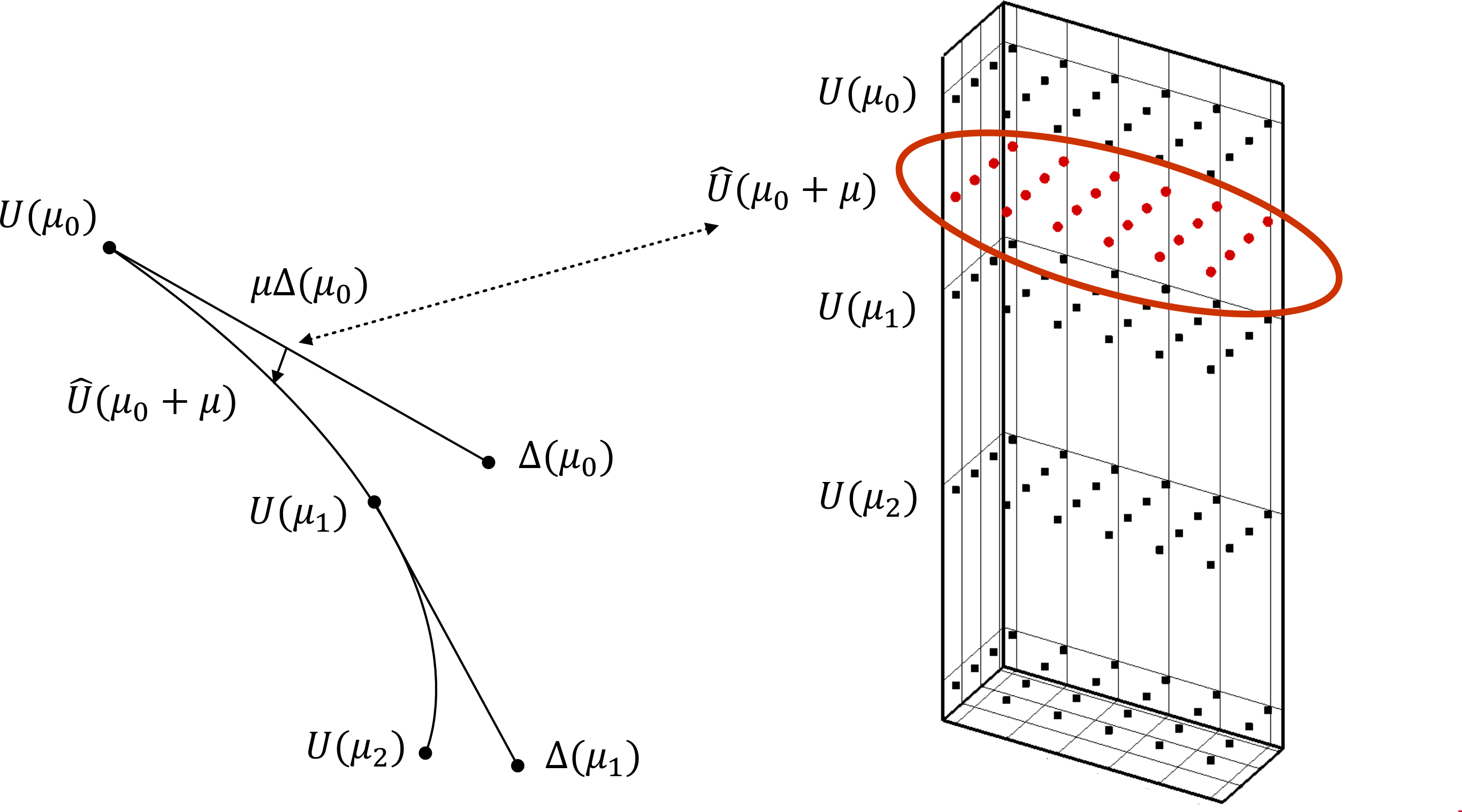}
\caption{Extrapolation of matrix manifold data. Sketched on the right is the sample matrix data in $\R^{n\times \p}$. The curved line on the left represents the nonlinear matrix manifold; the straight lines represent the tangent vectors in the tangent space.
The matrix curve is linearized at $U(q_0)$, $U(q_1)$, etc.}
\label{fig:ST_extra_scheme}
\end{figure}

{\em Performing the Taylor extrapolation on $St(n,\p)$.}
With $\U(\mu_0),\dot \U(\mu_0)$ at hand,
$\U(\mu_0+\mu)$ can be approximated using the Stiefel exponential:
$\U(\mu_0+\mu)\approx \hat{\U}(\mu_0+\mu):=Exp_{\U_0}^{St}(\mu \dot \U(\mu_0))$,
see Algorithm \ref{alg:Stexp}.
The process is illustrated in Fig. \ref{fig:ST_extra_scheme}.
%

Note that when the $\mu$-dependency is real-analytic, then the Euclidean Taylor expansion 
\begin{equation}
\label{eq:Taylor}
    \U(\mu_0 + \mu) =  \U(\mu_0) + \mu\dot{\U}(\mu_0)  + \frac{\mu^2}{2}\ddot{\U}(\mu_0)  + \mathcal{O}(\mu^3) \in St(n,\p)
\end{equation}
converges to an orthogonal matrix $\U(\mu_0 + \mu)\in St(n,\p)$.
Yet, when truncating the Taylor series, we leave the Stiefel manifold. 
In particular, the columns of the first order approximation are not orthonormal, i.e.
$\U(\mu_0) + \mu\dot{\U}(\mu_0) \notin St(n,\p)$ for $\mu \neq 0$.
By construction, the Stiefel geodesic 
features the same starting velocity $\dot{\U}(\mu_0)$ and thus matches the Taylor series up to terms of second order.
In addition, it respects the geometric structure of the Stiefel manifold and thus  preserves column-orthonormality for every $\mu$.
%
%
%
%
%
%
%
\section{Matrix manifolds of practical importance}
\label{sec:MatMan}
In this section, we discuss the matrix manifolds that feature most often in practical applications in the context of model reduction. For each manifold under consideration, we recap,  if applicable
\begin{itemize}
 \item the representation of points/locations in numerical schemes.
 \item the representation of tangent vectors in numerical schemes.
 \item the most common Riemannian metrics.
 \item how to compute distances, geodesics and the Riemannian exponential and logarithm mappings.
\end{itemize}
%
%
%
%
%
%
%
\subsection{The general linear group}
This section is devoted to the general linear group $GL(n)$ of invertible square matrices.
In model reduction, regular matrices appear for example as (reduced) system matrices in LTI and discretized PDE systems
\cite{AmsallemFarhat2011, Degroote_etal2010, Panzer_etal2010} and parameterizations have to be such that matrix regularity is preserved.
In addition, the discussion of the seemingly simple matrix manifold $GL(n)$ is important, because 
it is the fundamental matrix Lie Group from which all other matrix Lie groups are derived.
Moreover, it provides the background for understanding quotient spaces of $GL(n)$, see Subsection \ref{sec:LieGroups} and also \cite{BonnabelSepulchre2009, VandereyckenAV_2012}.
A short summary on the Riemannian geometry of $GL(n)$ is given in \cite[\S 6]{Rentmeesters2013}.
\subsubsection{Introduction and data representation in numerical schemes}
\label{sec:GL_intro_and_numerics}
%
%
Because $GL(n) = \det^{-1}(\R\setminus\{0\}) =  \{A\in \R^{n\times n}| \det(A)\neq 0\}$, $GL(n)$ is
an open subset of the $n^2$-dimensional vector space $\R^{n\times n}\simeq \R^{n^2}$ and is thus an $n^2$-dimensional differentiable manifold, see \cite[Examples 1.22--1.27]{Lee2012smooth}.
%
The matrix manifold $GL(n)$ is disconnected as it decomposes into two connected components, namely the regular matrices of positive determinant and the regular matrices of negative determinant.

Because $GL(n)$ is an open subset of the vector space $\R^{n\times n}$, the tangent space at a location $A\in GL(n)$ is simply
$T_AGL(n) = \R^{n\times n}$.
%
For $GL(n)$, \tc{the Lie algebra is $\mathfrak{gl}(n)=\R^{n\times n}$, so that} the Lie group exponential is the standard matrix exponential 
$\exp_m:\R^{n\times n}=\mathfrak{gl}(n)\rightarrow GL(n)$.
\tc{From the Lie group perspective \eqref{eq:tangspaceshifts}}, the tangent space at an arbitrary point 
$A\in GL(n)$ is to be considered as the set $T_AGL(n) = A\mathfrak{gl}(n)=A(\R^{n\times n})$, even though this set coincides with $\R^{n\times n}$.
\subsubsection{Distances and geodesics}
\label{sec:GL_data}
The obvious choice for a Riemannian metric on $GL(n)$ is to use the inner product from the ambient Euclidean matrix space, i.e.,
\[
 \langle \Delta, \tilde{\Delta}\rangle_A = \langle \Delta, \tilde{\Delta}\rangle_0 = \tr(\Delta^T \tilde{\Delta}),
\]
for $A\in GL(n)$ and $\Delta, \tilde{\Delta} \in T_AGL(n) =\R^{n\times n}$.

In many applications, it is more appropriate to consider metrics with certain invariance properties.\footnote{\tc{``Eulerian motion of a rigid body can be described as motion along geodesics
in the group of rotations of three-dimensional euclidean space provided with
a left-invariant Riemannian metric. A significant part of Euler's theory
depends only upon this invariance, and therefore can be extended to other
groups.''\cite[Appendix 2, p. 318]{Arnold1997mathematical}}}
A {\em left-invariant metric} can be obtained from the standard metric via
\begin{equation}
 \label{eq:GLn_leftinv_metric}
 \langle \Delta, \tilde{\Delta}\rangle_A =  \langle A^{-1}\Delta, A^{-1}\tilde{\Delta}\rangle_0, 
 \quad A\in GL(n), \quad \Delta, \tilde{\Delta} \in T_AGL(n).
\end{equation}
When formally considering $\Delta = AV, \tilde{\Delta}= A\tilde{V} \in T_AGL(n)=A\mathfrak{gl}(n)$ as left-translates of tangent vectors
$V,\tilde{V}\in T_IGL(n)=\mathfrak{gl}(n)$, 
then this metric satisfies
$\langle \Delta, \tilde{\Delta}\rangle_A =  \langle V, \tilde{V}\rangle_0$.
Alternatively, $\langle V, \tilde{V}\rangle_0 = \langle A V, A\tilde{V}\rangle_A$, which explains the name `left-invariant'.
\paragraph*{The Riemannian exponential and logarithm for the flat metric}
When equipped with the Euclidean metric, $GL(n)$ is flat: since the tangent space is the full matrix space $\R^{n\times n}$,
the geodesic equation \eqref{eq:GeodODE_euclid} requires the acceleration of a geodesic curve to vanish completely.
Hence, the geodesic that starts from $A\in GL(n)$ with velocity $\Delta \in \R^{n\times n}$ is the straight line
$C(t) = A + t\Delta$. Note that the curve $t\mapsto C(t)$ may leave the manifold $GL(n)$ for some $t\in\R$ as it may hit a matrix with zero determinant.
The formulae for the Riemannian exponential and logarithm mapping 
at a base point $A\in GL(n)$ are 
\begin{eqnarray}
 \label{eq:flat_exp_GLn}
 \Exp_A^{GL}:& T_AGL(n)\supset B_{\e}(0) \rightarrow GL(n),&
 \quad \Delta \mapsto \tilde{A} :=A+\Delta,\\
 \label{eq:flat_log_GLn}
 \Log_A^{GL}:& GL(n) \rightarrow  T_AGL(n),&
 \quad \tilde{A}\mapsto \Delta := (\tilde{A}-A).
\end{eqnarray}
In \eqref{eq:flat_exp_GLn}, $B_{\e}(0)$ denotes a suitably small open neighborhood around $0\in T_AGL(n)\simeq \R^{n\times n}$ such that $A + \Delta\in GL(n)$ for all $\Delta \in B_{\e}(0)$.
%
%
%
\paragraph*{The Riemannian exponential for the left-invariant metric on GL(n)}
The left-invariant metric induces a non-flat geometry on $GL(n)$.
Formulae for the covariant derivatives and the corresponding geodesics 
are derived in \cite[Thm. 2.14]{ANDRUCHOW2014241}.
The counterparts w.r.t. the right-invariant metrics can be found in 
\cite{VandereyckenAV_2012}.
Given a base point $A\in GL(n)$ and a starting velocity 
$\Delta = A V\in T_AGL(n) = A\mathfrak{gl}(n)$, the associated geodesic is 
\begin{equation}
 \label{eq:LI_Geo_GLn}
 \Gamma_{A,\Delta}: t\mapsto A \exp_m(tV^T)\exp_m(t(V-V^T)).
\end{equation}
The Riemannian exponential is
\begin{eqnarray}
\nonumber
 \Exp_M^{GL}(\Delta) &=& \Gamma_{A, \Delta}(1) = A \exp_m(V^T)\exp_m(V-V^T)\\
  \label{eq:Exp_LI_Geo_GLn}
 &=& A \exp_m((A^{-1}\Delta)^T)\exp_m((A^{-1}\Delta)-(A^{-1}\Delta)^T).
\end{eqnarray}
The author is not aware of a closed formula for the inverse map, i.e., the Riemannian logarithm for the left-invariant metric,
see also the discussion in \cite[\S 4.5]{VandereyckenAV_2012}.
The thesis \cite[\S 6.2]{Rentmeesters2013} introduces a Riemannian shooting method for computing the Riemannian logarithm w.r.t. the left-invariant metric.
\paragraph*{An important special case}
For tangent vectors $\Delta = AV \in T_AGL(n)$ with {\em normal} 
$V\in \R^{n\times n}$, i.e., $VV^T = V^TV$, it holds that the matrices
$V^T$ and $(V-V^T)$ commute.
Therefore, according to \eqref{eq:BCH}, 
$ A \exp_m(V^T)\exp_m(V-V^T) = A \exp_m(V^T+V-V^T)= A \exp_m(V)$ and
the Riemannian exponential reduces to
\[
 \Exp_A^{GL}:T_AGL(n)\cap \{\Delta| A^{-1}\Delta \mbox{ normal}\}\rightarrow GL(n),  \Delta\mapsto  \tilde{A} = A \exp_m(A^{-1}\Delta).
\]
The Riemannian logarithm is
\[
 \Log_A^{GL}: \mathcal{D}_A \cap \{\tilde{A}| A^{-1} \tilde{A} \mbox{ normal}\}\rightarrow T_AGL(n), 
 \quad \tilde{A} \mapsto \Delta =  A \log_m(A^{-1}\tilde{A}),
\]
where $\mathcal{D}_A\subset GL(n)$ is a domain such that a suitable branch of the matrix logarithm is well-defined.
These expressions are sometimes encountered in the literature as the Riemannian exponential and logarithm mappings.
Yet, one should be aware of the fact that they hold under special circumstances.
%
%
\subsection{The orthogonal group}
This section is devoted to the {\em orthogonal group} $O(n)\subset \R^{n\times n}$ of orthogonal $n$-by-$n$ matrices.
In parametric model reduction, such matrices may appear as eigenvector matrices in symmetric EVD problems.
\subsubsection{Introduction and data representation in numerical schemes}
\label{sec:On_intro_and_numerics}
The orthogonal group is $O(n)= \{Q\in\R^{n\times n}|\hspace{0.2cm} QQ^T=I=Q^TQ\}$.
The manifold structure of $O(n)$
can be established via Theorem \ref{thm:regurbild}, see also Example \ref{ex:reg_urbild_On}.
The orthogonal group decomposes into two connected components, namely the orthogonal matrices with 
determinant $1$ and the orthogonal matrices with determinant $-1$. The former constitute the 
{\em special orthogonal group} $SO(n) = \{Q\in O(n)| \det(Q)=1\}$.
The orthogonal group is a closed subgroup of the Lie group $GL(n)$ and thus itself a Lie group \tc{(Section \ref{sec:LieGroups})}.
%
%
The tangent space $T_IO(n)$ at the identity forms
the Lie algebra associated with the Lie group $O(n)$. It coincides with the Lie algebra of $SO(n)$
and as such is denoted by $\mathfrak{so}(n) = T_ISO(n) = T_IO(n)$, \cite[\S 3.3, 3.4]{Hall_Lie2015}.
The Lie algebra of $SO(n)$ is precisely the vector space of skew-symmetric matrices, $\mathfrak{so}(n)= \Skew(n)$.
\tc{According to \eqref{eq:tangspaceshifts},} the tangent space at an arbitrary location $Q$ is given by the translates (by left-multiplication) of the Lie algebra
\[
  T_QO(n) =Q \mathfrak{so}(n) = \left\{\Delta = QV\in \R^{n\times n}|\hspace{0.2cm} V\in \Skew(n)\right\},
\]
which is the same as $\left\{\Delta \in \R^{n\times n}|\hspace{0.2cm} Q^T\Delta = -\Delta^TQ\right\}$.
The Lie exponential is
\begin{equation}
 \label{eq:On_Lie_exp}
 \exp_m|_{\mathfrak{so}(n)}: \mathfrak{so}(n) \rightarrow SO(n).
\end{equation}
This restriction is a surjective map, see Appendix \ref{app:pre}.
The dimensions of both $T_QO(n)$ and $O(n)$ are $\frac{1}{2}n(n-1)$.
\subsubsection{Distances and geodesics}
\label{sec:On_data}
We follow up on the discussion in Section \ref{sec:GL_intro_and_numerics}.
For the orthogonal group, the Euclidean metric and the left-invariant metric coincide:
Let $\Delta = QV, \tilde{\Delta} = Q\tilde{V}\in T_QO(n) =Q \mathfrak{so}(n)$. Then,
\begin{eqnarray*}
 \langle \Delta, \tilde{\Delta} \rangle_Q =& \langle Q^{-1} \Delta, Q^{-1}\tilde{\Delta} \rangle_0
 &= \langle V, \tilde{V} \rangle_0\\
 =& \tr(V^T\tilde{V}) = \tr(V^TQ^TQ\tilde{V}) &= \langle \Delta, \tilde{\Delta} \rangle_I.
\end{eqnarray*}
In fact, the metric is also right-invariant, which makes it a {\em bi-invariant} metric, see \cite[\S 2]{Alexandrino2015}.
Bi-invariant metrics are important, because for Lie groups endowed with bi-invariant metrics,
the Lie exponential map and the Riemannian exponential map at the identity coincide \cite[Thm. 2.27, p. 40]{Alexandrino2015}.
\paragraph*{The Riemannian exponential and logarithm maps on O(n)}
The Riemannian $O(n)$-exponential at a base point $Q\in O(n)$ sends a tangent vector $\Delta\in T_QO(n)$
to the endpoint $\tilde{Q}\in O(n)$ of a geodesic that starts from $Q$ with velocity vector $\Delta$.
Therefore, it provides at the same time an expression for the geodesic curves on $O(n)$.
A formula for computing the Riemannian $O(n)$-exponential
was derived in \cite[\S 2.2.2]{EdelmanAriasSmith1999}.
Given $Q\in O(n)$, it holds
\begin{equation}
 \label{eq:On_exp}
 \Exp_Q^{On}: T_QO(n) \rightarrow O(n),\quad \Delta \mapsto \tilde{Q}:=Q\exp_m(Q^T\Delta).
\end{equation}
This result is also immediate from abstract Lie theory, see \cite[Eq. (2.2) \& Thm. 2.27]{Alexandrino2015}.\footnote{The Lie exponential is $\exp_m|_{\mathfrak{so}(n)}: \mathfrak{so}(n) \rightarrow SO(n)$, which is in the case at hand the Riemannian exponential at the identity, $\Exp_I^{SO} = \exp_m|_{\mathfrak{so}(n)}$. This translates to any other location via 
\cite[Eq. (2.2)]{Alexandrino2015} as follows: Pick any $Q\in SO(n)$ and consider the mapping ``left-multiplication by Q'', i.e., $L_Q: SO(n) \rightarrow SO(n), P\mapsto QP$. Then, the differential is $d(L_Q)_I: T_ISO(n) \rightarrow T_{L_Q(I)}SO(n), V\mapsto \Delta:= QV$. Because $L_Q$ is an isometry, 
\[
 Q\Exp_I^{SO}(V) = L_Q(\Exp_I^{SO}(V)) 
 = \Exp_{L_Q(I)}^{SO}(d(L_Q)_I(V)) = \Exp_Q^{SO}(QV),
\]
which gives $\Exp_{Q}^{SO}(QV) = Q\Exp_I^{SO}(V) = Q\exp_m(Q^{-1}\Delta)$ and thus \eqref{eq:On_exp}.
}
The corresponding Riemmanian logarithm on $O(n)$ is
\begin{equation}
 \label{eq:On_log}
 \Log_Q^{On}: O(n)\supset \mathcal{D}_Q \rightarrow T_QO(n),\quad \tilde{Q} \mapsto \Delta:= Q\log_m(Q^T\tilde{Q})
\end{equation}
and is well defined on a neighborhood $\mathcal{D}_Q\subset O(n)$ around $Q$ such that for all $\tilde{Q}\in \mathcal{D}_p$, the orthogonal matrix $Q^T\tilde{Q}$ does not feature $\lambda = -1$ as an eigenvalue.
%
%
%
\paragraph*{The Riemannian distance between orthogonal matrices}
For given $Q, \tilde{Q}\in O(n)$ from the same connected component of $O(n)$,
consider the EVD $Q^T\tilde{Q} = \Psi \Lambda \Psi^H$.
Because $Q^T\tilde{Q}$ is orthogonal, it holds $\Lambda = \diag(e^{i\theta_1}, \ldots, e^{i\theta_n})$ 
and we assume that $\theta_1,\ldots,\theta_n\in (-\pi, \pi)$.
The Riemannian distance is
\begin{align*}
 \label{eq:RiemannDist_On}
 \dist_{On}(Q, \tilde{Q}) &= \| \Log_Q^{On}(\tilde{Q})\|_Q 
 = \|\log_m (\Lambda)\|_F
 = \left(\sum_{k=1}^n \theta_k^2\right)^{\frac{1}{2}}.
\end{align*}
The compact Lie group $SO(n)$ is a geodesically complete Riemannian manifold \cite[Hopf-Rinow-Theorem, p. 31]{Alexandrino2015}, and each two points of $SO(n)$ can be joined by a minimal geodesic.
%
%
%
%
%
\subsection{The matrix manifold of symmetric positive definite matrices}
%
%
%
This section is devoted to the matrix manifold $SPD(n)$ of real, symmetric positive-definite $n$-by-$n$ matrices.
In model reduction, such matrices appear for example as (reduced) system matrices in second-order parametric ODEs.
For example, in linear structural or electrical dynamical systems, mass, stiffness and damping matrices are usually 
in $SPD(n)$, \cite[\S 4.2]{AmsallemFarhat2011}. 
\tc{Moreover, positive definite matrices arise as Gramians of reachable and observable LTI systems in the context of balanced truncation \cite{BennerGugercinWillcox2015}.} 
Related is the manifold of positive {\em semi-definite} matrices of fixed rank.
It is investigated in \cite{BonnabelSepulchre2009, VandereyckenAV_2012, Massart2018b}.
An application in the context of model reduction features in \cite{Massart2018}.
\subsubsection{Introduction and data representation in numerical schemes}
\label{sec:SPD_intro_and_numerics}
The set 
\[
 SPD(n) = \{A\in\sym(n)|\hspace{0.2cm} x^TAx > 0 \quad \forall x\in\R^n\setminus\{0\}\}
\]
is an open subset of the metric Hilbert space $(\sym(n), \langle\cdot,\cdot\rangle_0)$ of symmetric matrices.
As such, it is a differentiable manifold \cite[\S 6]{Bhatia2007}.
Moreover, it forms a {\em convex cone} \cite[Example 2, p. 8]{Faraut1994}, \cite[\S 2.3]{Moakher2005}, and can be realized as a quotient $SPD(n) \simeq GL(n)/O(n)$.
The latter is based on the fact that for $A\in SPD(n)$, matrix factorizations $A = ZZ^T$ with $Z\in GL(n)$ 
are invariant under orthogonal transformations $Z\mapsto ZQ$, $Q\in O(n)$, \cite[\S 2, p.3]{BonnabelSepulchre2009}.

Since $SPD(n)$ is an open subset of the vector space $\sym(n)$, the tangent space is simply
\begin{equation}
 \label{eq:spd_tnagspace}
 T_ASPD(n) = \sym(n).
\end{equation}
The dimensions of both $T_ASPD(n)$ and $SPD(n)$ are $\frac{1}{2}n(n+1)$.

There is a smooth one-to-one correspondence between $\sym(n)$ and $SPD(n)$.
That is, every positive definite matrix can be written as the matrix exponential of a unique symmetric matrix,
\cite[Lem. 18.7, p. 472]{Gallier2011}.
Put in different words, when restricted to $\sym(n)$, the standard matrix exponential
\[
 \exp_m: \sym(n) \rightarrow SPD(n)
\]
is a diffeomorphism, its inverse is the standard principal matrix logarithm
\[
 \log_m: SPD(n)\rightarrow \sym(n),
\]
see also \cite[Thm. 2.8]{ArsignyFPA06}.
The group $GL(n)$ acts on $SPD(n)$ via congruence transformations
\begin{equation}
\label{eq:spd_cong}
  g_X(A) = X^TAX, \quad X\in GL(n), A\in SPD(n).
\end{equation}
For additional background on $SPD(n)$, see \cite{Moakher2006, Moakher2011, Pennec2006}.
Applications in computer vision are presented in \cite{CherianSra_2016, Kim_et_al_CCA_on_SPDn2015}.
\subsubsection{Distances and geodesics}
\label{sec:SPD_dist}
The literature knows a large variety of distance measures on $SPD(n)$, see \cite[Table 3.1, p. 56]{RiemannKernels_Jayasumana2015}.
Yet, there are essentially two choices that are associated with inner products on the tangent space of $SPD(n)$ and thus 
induce Riemannian geometries on the manifold $SPD(n)$:
the so-called \tc{\em natural metric} and the {\em log-Euclidean metric}.
Let $A\in SPD(n)$ and let $\Delta,\tilde{\Delta}\in \sym(n)$ be two tangent vectors.
\begin{itemize}
 \item The \tc{natural metric} is
  \[
   \langle \Delta, \tilde{\Delta} \rangle_A 
   = \langle A^{-1/2}\Delta A^{-1/2}, A^{-1/2}\tilde{\Delta}A^{-1/2} \rangle_0
   = \tr(A^{-1}\Delta A^{-1}\tilde{\Delta}),
  \]
  see \cite[\S 6, p. 201]{Bhatia2007}, \cite{BonnabelSepulchre2009}. 
  It also goes by the name {\em trace matric}, \cite[\S XII.1, p.322]{lang2001fundamentals}.
  In statistical applications, it is usually called the {\em affine-invariant metric} \cite{MinhMurino2016, Pennec2006b}.\footnote{\tc{
  The motivation is as follows: if $y = Ax + v_0$, $A\in GL(n)$ is an affine transformation of a random vector $x$,
  then the mean is transformed to $\bar{y} := A\bar{x} + v_0$ and the covariance matrix undergoes a congruence transformation  $C_{yy} = E[(y-\bar{y})(y-\bar{y})^T ] = A C_{xx}A^T$.}}
   
 \item The log-Euclidean metric is
   \[
   \langle \Delta, \tilde{\Delta} \rangle_A = \langle D(\log_m)_A(\Delta), D(\log_m)_A(\tilde{\Delta})\rangle_0,
  \]
  see \cite[eq. (3.5)]{ArsignyFPA06}.
\end{itemize}
  For the natural metric, it is more appropriate to consider 
  $\sym(n) = T_ISPD(n)$ as the tangent space at the identity and the tangent space at an arbitrary
  location $A\in SPD(n)$ as $T_ASPD(n) = A^{1/2}\left(T_ISPD(n)\right) A^{1/2}$,
  which, of course, is nothing but a reparameterization of $\sym(n)$.
  From this perspective, we have for tangent vectors
  $\Delta = A^{1/2}VA^{1/2}, \tilde{\Delta}= A^{1/2}\tilde{V}A^{1/2}$
  that 
  \[
   \langle \Delta, \tilde{\Delta} \rangle_A = \langle V, \tilde{V} \rangle_0.
  \]

\tc{The congruence transformations \eqref{eq:spd_cong} are isometries of $SPD(n)$ with respect to the natural metric}, \cite[Thm. XII.1.1, p. 324]{lang2001fundamentals}, \cite[Lem. 6.1.1, p. 201]{Bhatia2007}.
See also the discussion in \cite[\S 3]{Pennec2006b}.

By a standard pullback construction from differential geometry \cite[Def. 2.2, Example 2.5]{DoCarmo2013riemannian}, the log-Euclidean metric transfers the inner product 
$\langle \cdot,\cdot\rangle_0$ on $\sym(n)$ to $SPD(n)$ via the matrix logarithm
$\log_m: SPD(n)\rightarrow \sym(n)$.
In \cite[eq. (3.5)]{ArsignyFPA06}, the authors take this construction one step further 
and use the $\exp_m$-$\log_m$-correspondence to define a multiplication that turns
$SPD(n)$ into a Lie group and, eventually, into a vector space.
As such, it is a {\em flat} manifold, i.e. a Riemannian manifold with zero curvature.
In this way, the computational challenges that come with dealing with data on nonlinear manifolds 
are circumvented.

Which metric is to be preferred is problem-dependent, see the various contributions in 
\cite{RiemannInComputerVision} and \cite{Minh:2016:AAR:3029338}.
Since the natural metric arises canonical both from the geometric approach, \cite[\S XII.1]{lang2001fundamentals},
and the matrix-algebraic approach \cite[\S 6]{Bhatia2007} and 
since staying with the standard matrix multiplication is consistent with 
the setting of solving dynamical systems in model reduction applications,
we restrict the discussion of the Riemannian exponential and logarithm to the geometry that is based on the natural metric. 
\paragraph*{The SPD(n) exponential}
The Riemannian $SPD(n)$-exponential at a base point $A\in SPD(n)$ sends a tangent vector $\Delta$
to the endpoint $\tilde{A}\in SPD(n)$ of a geodesic that starts from $A$ with velocity vector $\Delta$.
Therefore, it provides at the same time an expression for the geodesic curves on $SPD(n)$ with respect to the natural metric.
Formulae for computing the $SPD(n)$-exponential can be found in \cite{BonnabelSepulchre2009}, \cite{Pennec2006b}.
Readers preferring a matrix-analytic approach are referred to \cite[\S 6]{Bhatia2007}.
\begin{algorithm}
\caption{Riemanian $SPD(n)$-exponential}
\label{alg:SPDexp}
\begin{algorithmic}[1]
  \REQUIRE{base point $A\in SPD(n)$,  tangent vector $\Delta\in T_ASPD(n)=\sym(n)$}
  \ENSURE{$ \tilde{A}:= \Exp_{A}^{SPD}(\Delta) = A^{\frac12} \exp_m\left(A^{-\frac12} \Delta A^{-\frac12}\right)A^{\frac12}$.}
\end{algorithmic}
\end{algorithm}
Here,  $A^{\frac12}$ denotes the matrix square root of $A$, see Appendix \ref{app:pre}.
\paragraph*{The SPD(n) logarithm}
The Riemannian $SPD(n)$-logarithm at a base point $A\in SPD(n)$ finds for another point $\tilde{A}\in SPD(n)$ 
an $SPD(n)$-tangent vector $\Delta$ such that the geodesic that starts from $A$ with velocity $\Delta$ reaches 
$\tilde{A}$ after an arc length of $\|\Delta\|_A =\sqrt{\langle \Delta,\Delta \rangle_A}$.
Therefore, it provides for two given data points $A,\tilde{A}\in SPD(n)$
\begin{itemize}
 \item a solution to the geodesic endpoint problem: a geodesic that starts from $A$ and ends at $\tilde{A}$.
 \item the Riemannian distance between the given points $A, \tilde{A}$.
\end{itemize}
Formulae for computing the $SPD(n)$-logarithm can be found in \cite{BonnabelSepulchre2009}, \cite{Pennec2006b}.
\begin{algorithm}
\caption{Riemanian $SPD(n)$-logarithm}
\label{alg:SPDlog}
\begin{algorithmic}[1]
  \REQUIRE{base point $A\in SPD(n)$,  location $\tilde{A}\in SPD(n)$}
  \ENSURE{$ \Delta:= \Log_{A}^{SPD}(\tilde{A}) = A^{\frac12} \log_m\left(A^{-\frac12} \tilde{A} A^{-\frac12}\right)A^{\frac12}$.}
\end{algorithmic}
\end{algorithm}
Both Algorithms \ref{alg:SPDexp} and \ref{alg:SPDlog} require to compute the spectral decomposition of $n$-by-$n$-matrices. The computational effort is $\mathcal{O}(n^3)$.
In the context of parametric model reduction, the Riemannian exponential and logarithm maps are usually required for reduced matrix operators
\cite{AmsallemFarhat2011}. If $n$ denotes the dimension of the full state vectors and
$\p\ll n$ denotes the dimension of the reduced state vectors, then matrix exponentials for $\p$-by-$\p$-matrices are required, so that the computational effort reduces to $\mathcal{O}(\p^3)$.
%
%
%
%
%
%
%
\subsection{The Stiefel manifold}
\label{sec:Stiefel}
%
%
%
This section is devoted to the {\em Stiefel manifold} $St(n,\p)\subset \R^{n\times \p}$ of rectangular column-orthogonal $n$-by-$\p$ matrices, $\p\leq n$.
Points $U\in St(n,\p)$ may be considered as orthonormal bases of cardinality $\p$,  or $\p$-frames in $\R^n$.
In model reduction, such matrices appear as orthogonal coordinate systems for low-order ansatz spaces that usually stem from a proper orthogonal decomposition or a singular value decomposition of given input solution data.
Modeling data on the Stiefel manifold corresponds to data processing for orthonormal bases 
and thus allows for example for interpolation/parameterization of POD subspace bases.
The most important use case in model reduction is where the Stiefel matrices are tall and skinny, i.e., $\p\ll n$.
Interpolation problems on the Stiefel manifold have not yet been considered in the model reduction context.
The reference \cite{Krakowski2015SOLVINGIP} discusses interpolation of Stiefel data, however with using quasi-geodesics rather than geodesics. 
The work \cite{Zimmermann_2019} includes numerical experiments for interpolating orthogonal frames on the Stiefel manifold that 
relies the canonical Riemannian Stiefel logarithm \cite{Rentmeesters2013, StiefelLog_Zimmermann2017}.
\subsubsection{Introduction and data representation in numerical schemes}
\label{sec:Stiefel_intro_and_numerics}
The {\em Stiefel manifold} is the compact, homogeneous matrix manifold of column-orthogonal
matrices
\[
  St(n,\p):= \{U \in \R^{n\times \p}| \hspace{0.1cm}  U^TU = I_\p\}.
\]
The manifold structure can be directly established via Theorem \ref{thm:regurbild} in a similar way as in Example \ref{ex:reg_urbild_On}.
\tc{An alternative approach is via Example \ref{ex:QuotBasic}, where $St(n,\p)$ is identified with the quotient space
$St(n,\p)\cong O(n)/(I_\p\times O(n-\p))$
under actions of the closed subgroup $I_\p\times O(n-\p):= \left\{\begin{pmatrix}
                            I_\p &0\\
                             0    & R
                           \end{pmatrix}|\hspace{0.1cm} R\in O(n-\p)\right\}\leq O(n)$.
Two square orthogonal matrices in $O(n)$ are identified 
as the same point on $St(n,\p)$, if their first $\p$ columns coincide, see \cite[\S 2.4]{EdelmanAriasSmith1999}.}

For any matrix representative $U\in St(n,\p)$,
the tangent space of $St(n,\p)$ at $U$ is represented by
\[
  T_USt(n,\p) = \left\{\Delta \in \R^{n\times \p}|\hspace{0.2cm} U^T\Delta = -\Delta^TU\right\}\subset \R^{n\times \p}.
\]
Every tangent vector $\Delta \in T_USt(n,\p)$ may be written as
\begin{eqnarray}
\label{eq:tang_St}
  \Delta &=& UA + (I-UU^T)T, \quad A \in \R^{\p\times \p} \mbox{ skew}, \quad T\in\R^{n\times \p} \mbox{arbitrary,}\\
  \Delta &=& UA + U^{\bot}B, \quad A \in \R^{\p\times \p} \mbox{ skew}, \quad B\in\R^{(n-\p)\times \p} \mbox{ arbitrary,}
\end{eqnarray}
where in the latter case, $U^{\bot}\in St(n,n-\p)$ is such that $(U, U^{\bot})\in O(n)$ is a square orthogonal matrix.
The dimension of both $T_USt(n,\p)$ and $St(n,\p)$ is $n\p -\frac{1}{2}\p(\p+1)$.
For additional background and applications, see \cite{AbsilMahonySepulchre2008, bernstein2012tangent, Chakraborty2017, EdelmanAriasSmith1999, HueperRollingStiefel2008, Turaga_2008}.
\subsubsection{Distances and geodesics}
\label{sec:Stiefel_data}
Let $U\in St(n,\p)$ be a point and let $\Delta= UA + U^{\bot}B$, $\tilde{\Delta}= U\tilde{A} + U^{\bot}\tilde{B} \in T_USt(n,\p)$ be tangent vectors.
There are two standard metrics on the Stiefel manifold.
\begin{itemize}
 \item The {\em Euclidean metric} on $T_USt(n,\p)$ is the one inherited from the ambient $\mathbb{R}^{n\times \p}$:
 \[
   \langle\Delta, \tilde{\Delta}\rangle_0 = \tr( \Delta^T \tilde{\Delta}) = \tr A^T\tilde{A} + \tr B^T\tilde{B}
 \]
 \item The {\em canonical metric} on $T_USt(n,\p)$ 
  \[
  \langle \Delta, \tilde{\Delta}\rangle_U = \tr\left(\Delta^T(I-\frac{1}{2}UU^T)\tilde{\Delta}\right) 
  = \frac12 \tr A^T\tilde{A} +  \tr B^T\tilde{B}
 \]
 is derived from the
quotient representation $St(n,\p) = O(n)/(I_\p\times O(n-\p))$ of the Stiefel manifold.
\end{itemize}
The canonical metric counts the {\em independent} coordinates\footnote{i.e., the upper triangular entries of the skew-symmetric $A$ and the entries of $B$ of $\Delta = UA+U^{\bot}B$} of a tangent vector equally,
when measuring the length $\sqrt{\langle \Delta,\Delta\rangle_U}$ of a tangent vector
$\Delta = UA+U^{\bot}B$, while the Euclidean metric disregards the 
skew-symmetry of $A$ \cite[\S 2.4]{EdelmanAriasSmith1999}.
Recall that different metrics entail different measures for the lengths of curves and thus different formulae for geodesics.
\paragraph*{The Stiefel exponential}
The Riemannian Stiefel exponential at a base point $U\in St(n,\p)$ sends a Stiefel tangent vector $\Delta$
to the endpoint $\tilde{U}\in St(n,\p)$ of a geodesic that starts from $U$ with velocity vector $\Delta$.
Therefore, it provides at the same time an expression for geodesic curves on $St(n,\p)$.

\tc{A closed-form expression for the Stiefel exponential w.r.t. Euclidean metric is included in 
\cite[\S 2.2.2]{EdelmanAriasSmith1999},
  \[
   \tilde{U} = \Exp_U^{St}(\Delta) = 
   \left(U, \Delta\right)
   \exp_m\left(
     \begin{pmatrix}
   U^T\Delta & -\Delta^T\Delta\\
   I_p & U^T\Delta
  \end{pmatrix}
   \right)
     \begin{pmatrix}
   I_p\\
   0
  \end{pmatrix}
  \exp_m(-U^T\Delta).
  \]
In \cite{HueperUllrich2018}, an alternative formula is derived that features only matrix exponentials of skew-symmetric matrices.}
An efficient algorithm for computing the Stiefel exponential w.r.t. the canonical metric
was derived in \cite[\S 2.4.2]{EdelmanAriasSmith1999}:
\begin{algorithm}
\caption{Stiefel exponential \cite{EdelmanAriasSmith1999}.}
\label{alg:Stexp}
\begin{algorithmic}[1]
  \REQUIRE{base point $U\in St(n,\p)$,  tangent vector $\Delta\in T_USt(n,\p)$}
  \STATE{ $A := U^T\Delta$}  \hfill \COMMENT{horizontal component, skew}
  \STATE{ $QR := \Delta-UA$}  \hfill \COMMENT{(thin) qr-decomp. of normal component of $\Delta$.}
    \STATE{$\left(
          \begin{array}{cc}
            A & -R^T \\
            R & 0
          \end{array}
          \right) = T \Lambda T^H \in \R^{2\p \times 2\p}$}
   \hfill \COMMENT{EVD of skew-symmetric matrix}
  \STATE{ $\begin{pmatrix}M\\ N\end{pmatrix}
          := T\exp_m(\Lambda)T^H 
          \begin{pmatrix} I_\p\\\mathbf{0}\end{pmatrix}\in \R^{2\p\times \p}$ }
  \ENSURE{$\tilde{U} := Exp_{U}^{St}(\Delta) = UM + QN \in St(n,\p)$}
\end{algorithmic}
\end{algorithm}
%
In applications, where $Exp_{U}^{St}(\mu\Delta)$ needs to be evaluated for various parameters $\mu$ as in
in the example of Section \ref{sec:extrapol}, 
steps 1.--3. should be computed a priori (offline).
Apart from elementary matrix multiplications, the algorithm requires to compute the standard matrix exponential of a skew-symmetric matrix. This however, is for a $2\p$-by-$2\p$-matrix and does not scale in the dimension $n$.
With the usual assumption of model reduction that $n\gg p$, the computational effort is $\mathcal{O}(n\p^2)$.
\paragraph*{The Stiefel logarithm}
The Riemannian Stiefel logarithm at a base point $U\in St(n,\p)$ finds for another point $\tilde{U}\in St(n,\p)$ 
a Stiefel tangent vector $\Delta$ such that the geodesic that starts from $U$ with velocity $\Delta$ reaches 
$\tilde{U}$ after an arc length of $\|\Delta\|_U =\sqrt{\langle\Delta,\Delta\rangle_U}$.
Therefore, it provides for two given data points $U,\tilde{U}\in St(n,\p)$
\begin{itemize}
 \item a solution to the geodesic endpoint problem: a geodesic that starts from $U$ and ends at $\tilde{U}$.
 \item the Riemannian distance between the given points $U, \tilde{U}$.
\end{itemize}

An efficient algorithm for computing the Stiefel logarithm w.r.t. the canonical metric
was derived in \cite{StiefelLog_Zimmermann2017}.
\begin{algorithm}
\caption{Stiefel logarithm \cite{StiefelLog_Zimmermann2017}.}
\label{alg:Stlog}
\begin{algorithmic}[1]
  \REQUIRE{base point $U\in St(n,\p)$,  $\tilde{U}\in St(n,\p)$ `close' to base point, $\tau>0$
  convergence threshold}
  \STATE{ $M:= U^T\tilde{U} \in \R^{\p\times \p}$}
  \STATE{ $QN := \tilde{U} - UM \in \R^{n \times \p}$}
    \hfill \COMMENT{(thin) qr-decomp. of normal component of $\tilde{U}$}
  \STATE{ $V_0 := \begin{pmatrix}M&X_0\\N&Y_0\end{pmatrix} \in O(2\p)$}
    \hfill \COMMENT{compute orth. completion of the block $\begin{pmatrix}M\\N\end{pmatrix}$}

  \FOR{ $k=0,1,2,\ldots$}
  \STATE{$\begin{pmatrix}A_k & -B_k^T\\ B_k & C_k\end{pmatrix} := \log_m(V_k)$}
    \hfill \COMMENT{matrix log of orth. matrix}
  \IF{$\|C_k\|_2 \leq \tau$}
    \STATE{break}
  \ENDIF
  \STATE{ $\Phi_{k} := \exp_m{(-C_k)}$}
    \hfill \COMMENT{matrix exp of skew matrix}
  \STATE{$V_{k+1} := V_{k}W_k$, where $W_k:=\begin{pmatrix}I_\p& 0\\ 0 & \Phi_k\end{pmatrix}$}
  \ENDFOR
  \ENSURE{$\Delta := Log_{U}^{St}(\tilde{U}) = U A_k + QB_k \in T_{U}St(n,\p)$}
\end{algorithmic}
\end{algorithm}
The analysis in \cite{StiefelLog_Zimmermann2017} shows that the algorithm is guaranteed to converge if
the input data points $U,\tilde{U}$ are at most a Euclidean distance of $d=\|U-\tilde{U}\|_2\leq 0.09$ apart.
In this case, the algorithm exhibits a linear rate of convergence that depends on $d$ but is smaller than $\frac12$.
In practice, the algorithm seems to converge, whenever the initial $V_0$ is such that its standard matrix logarithm
$\log_m(V_0)$ is well-defined. Note that two points on $St(n,\p)$ can at most be a Euclidean distance of $2$ away from each other.

Apart from elementary matrix multiplications, the algorithm requires to compute the
standard matrix logarithm of an orthogonal $2\p$-by-$2\p$-matrix and the standard matrix exponential of a skew-symmetric
$\p$-by-$\p$-matrix at every iteration $k$. Yet, these operations are independent of the dimension $n$.
With the usual assumption of model reduction that $\p\ll n$, the computational effort is $\mathcal{O}(n\p^2)$.

For the Stiefel manifold equipped with the Euclidean metric, methods for calculating the Stiefel logarithm are introduced in
\cite{Bryner2017}.
%
%
%
%
%
%
%
\subsection{The Grassmann manifold}
\label{sec:Grassmann}
%
%
%
This section is devoted to the {\em Grassmann manifold} $Gr(n,\p)$ of $\p$-dimensional subspaces of $\R^n$ for $\p\leq n$.
Every point $\mathcal{U}\in Gr(n,\p)$, i.e., every subspace may be represented by selecting a basis
$\{u^1,\ldots,u^\p\}$ with $\colspan(u^1,\ldots,u^\p) = \mathcal{U}$.
In numerical schemes, we work exclusively with orthonormal bases. In this way, points $\mathcal{U}$ on the Grassmann manifold are to be represented by points $U\in St(n,\p)$ on the Stiefel manifold via $\mathcal{U} = \colspan(U)$.
For details and theoretical background, see the references 
\cite{AbsilMahonySepulchre2004, AbsilMahonySepulchre2008, EdelmanAriasSmith1999}.
Subspaces and Grassmann manifolds play an important role in projection-based parametric model reduction, \cite{AmsallemFarhat2008,Nguyen2012,Zimmermann2014,TanhSon2019}
and in Krylov subspace approaches \cite{BennerGugercinWillcox2015}. 
Modeling data on the Grassmann manifold corresponds to data processing for subspaces 
and thus allows for example for the interpolation/parameterization of POD subspaces.
The most important use case in model reduction is where the subspaces are of low dimension when compared to the 
surrounding state space, i.e., $n\gg p$.

\subsubsection{Introduction and data representation in numerical schemes}
\label{sec:Grassmann_intro_and_numerics}
The set of all $\p$-dimensional subspaces $\mathcal{U}\subset \R^n$
forms the {\em Grassmann manifold}
\[
  Gr(n,\p):= \{\mathcal{U}\subset \R^n| \hspace{0.1cm}  \mathcal{U} \mbox{ subspace, dim}(\mathcal{U}) = \p\}.
\]
\tc{
The Grassmann manifold is a quotient of $O(n)$ under the action of the Lie subgroup
$O(\p)\times O(n-\p) = 
\left\{
\begin{pmatrix}
 S & 0\\
 0 & R
\end{pmatrix}|\hspace{0.1cm} S\in O(\p), R\in O(n-\p)
\right\}\leq O(n)$.
Two matrices $Q, \tilde{Q}\in O(n)$ are in the same $(O(\p)\times O(n-\p))$-orbit, if and only if 
the first $\p$ columns of $Q$ and $\tilde{Q}$ span the same subspace and the tailing $n-\p$ columns 
span the corresponding orthogonal complement subspace. Theorem \ref{thm:HomSpaceConst} applies and shows that 
$Gr(n,\p)= O(n)/(O(\p)\times O(n-\p))$ is a homogeneous manifold.
}

Alternatively, the Grassmann manifold can be realized as a quotient manifold
of the Stiefel manifold \tc{with the help of Theorem \ref{thm:QuotMan}},
\begin{equation}
\label{eq:Grassmann_quotient}
  Gr(n,\p) = St(n,\p)/O(\p) = \{[U]| \hspace{0.1cm}  U\in St(n,\p)\},
\end{equation}
where the $O(\p)$-orbits are $[U] = \{UR| \hspace{0.1cm}  R\in O(\p)\}$.
A matrix $U\in St(n,\p)$ is called a {\em matrix representative} of a subspace $\mathcal{U}\in Gr(n,\p)$,
if $\mathcal{U}=\colspan(U)$.
The orbit $[U]$ and the subspace $\mathcal{U}=\colspan(U)$ are to be considered as the same object.
For any matrix representative $U\in St(n,\p)$ of $\mathcal{U}\in Gr(n,\p)$
the tangent space of $Gr(n,\p)$ at $\mathcal{U}$ is represented by
\[
 T_{\mathcal{U}}Gr(n,\p) = \left\{\Delta \in \R^{n\times \p}|\hspace{0.2cm} U^T\Delta = 0\right\}\subset \R^{n\times \p}.
\]
Every tangent vector $\Delta \in T_\mathcal{U}Gr(n,\p)$ may be written as
\begin{eqnarray}
\label{eq:tang_Gr}
  \Delta &=& (I-UU^T)T, \quad T\in\R^{n\times \p} \mbox{ arbitrary, or, }\\
  \Delta &=& U^{\bot}B, \quad B\in\R^{(n-\p)\times \p} \mbox{ arbitrary,}
\end{eqnarray}
where in the latter case, $U^{\bot}\in St(n,n-\p)$ is such that $(U, U^{\bot})\in O(n)$ is a square orthogonal matrix.
The dimension of both $T_\mathcal{U}Gr(n,\p)$ and $Gr(n,\p)$ is $n\p-\p^2$.
\subsubsection{Distances and geodesics}
\label{sec:Grassmann_dist}
A metric on $T_\mathcal{U}Gr(n,\p)$ can be obtained via making use of the fact that the Grassmannian is a quotient of the Stiefel manifold. Alternatively, one can restrict the standard inner matrix product $\langle A, B\rangle_0 = \tr(A^TB)$ to the Grassmann tangent space.
In the case of the Grassmannian, both approaches lead to the same metric
\[
 \langle \Delta, \tilde{\Delta}\rangle_\mathcal{U} =  \tr(\Delta^T\tilde{\Delta})= \langle \Delta, \tilde{\Delta}\rangle_0 ,
\]
see \cite[\S 2.5]{EdelmanAriasSmith1999}.
\paragraph*{The Grassmann exponential}
The Riemannian Grassmann exponential at a base point $\mathcal{U}\in Gr(n,\p)$ sends a Grassmann tangent vector $\Delta$
to the endpoint $\tilde{\mathcal{U}}\in Gr(n,\p)$ of a geodesic that starts from $\mathcal{U}$ with velocity vector $\Delta$.
Therefore, it provides at the same time an expression for the geodesic curves on $Gr(n,\p)$.
An efficient algorithm for computing the Grassmann exponential
was derived in \cite[\S 2.5.1]{EdelmanAriasSmith1999}:
\begin{algorithm}
\caption{Grassmann exponential \cite{EdelmanAriasSmith1999}.}
\label{alg:Grexp}
\begin{algorithmic}[1]
  \REQUIRE{base point $\mathcal{U}=[U]\in Gr(n,\p)$, where $U\in St(n,\p)$,  tangent vector $\Delta\in T_UGr(n,\p)$}
  \STATE{ $Q \Sigma V^T \stackrel{\mbox{\footnotesize SVD}}{:=} \Delta$, with $Q\in St(n,\p)$}  
  \hfill \COMMENT{(thin) SVD of tangent vector}
  \STATE{ $ \tilde{U} :=  UV\cos(\Sigma)V^T + Q\sin(\Sigma)V^T$}  
  \hfill \COMMENT{$\cos$ and $\sin$ act only on diag. entries.}
  \ENSURE{$ \tilde{\mathcal{U}}:= Exp_{\mathcal{U}}^{Gr}(\Delta) = [\tilde{U}]\in Gr(n,\p)$.}
\end{algorithmic}
\end{algorithm}
Apart from elementary matrix multiplications, the algorithm requires to compute the singular value decomposition of an
$n$-by-$\p$-matrix. The computational effort is $\mathcal{O}(n\p^2)$.
\paragraph*{The Grassmann logarithm}
The Riemannian Grassmann logarithm at a base point $\mathcal{U}\in Gr(n,\p)$ finds for another point $\tilde{\mathcal{U}}\in Gr(n,\p)$ 
a Grassmann tangent vector $\Delta$ such that the geodesic that starts from $\mathcal{U}$ with velocity $\Delta$ reaches 
$\tilde{\mathcal{U}}$ after an arc length of $\|\Delta\|_\mathcal{U} =\sqrt{g^C_\mathcal{U}(\Delta,\Delta)}$.
Therefore, it provides for two given data points $\mathcal{U},\tilde{\mathcal{U}}\in Gr(n,\p)$
\begin{itemize}
 \item a solution to the geodesic endpoint problem:  a geodesic that starts from $\mathcal{U}$ and ends at $\tilde{\mathcal{U}}$.
 \item the Riemannian distance between the given points $\mathcal{U}, \tilde{\mathcal{U}}$.
\end{itemize}
An algorithm for computing the Grassmann logarithm is stated implicitly in \cite[\S 3.8, p. 210]{AbsilMahonySepulchre2004}.
The reference \cite{Gallivan_etal2003} features expressions for the Grassmann
exponential and the corresponding logarithm that formally work with Grassmann representatives in $SO(n)/(SO(\p)\times SO(n-\p))$ but also keep the computational effort $\mathcal{O}(n\p^2)$.
The reference \cite[\S 4.3]{Rahman_etal2005} gives the corresponding mappings after identifying
subspaces with orthoprojectors, see also \cite{Batzies2015}.
\begin{algorithm}
\caption{Grassmann Logarithm.}
  \label{alg:Glog}
  \begin{algorithmic}[1]
  \REQUIRE{base point $\mathcal{U}=[U]\in G(n,\p)$ with $U\in St(n,\p)$, 
           $\tilde{\mathcal{U}}= [\tilde{U}]\in G(n,\p)$ with $\tilde{U}\in St(n,\p)$.}
  \STATE{ $M:= U^T\tilde{U}$}
  \STATE{ $L := (I -UU^T) \tilde{U} M^{-1} = \tilde{U} M^{-1} - U$ }
  \STATE{ $Q\Sigma V^T \stackrel{SVD}{:=} L$}
  \hfill \COMMENT{(thin) SVD }
  \STATE{ $\Delta := Q\arctan(\Sigma)V^T$}
  \hfill \COMMENT{$\arctan$ acts only on diag. entries.}
  \ENSURE{$\Delta = \Log_\mathcal{U}^{Gr}(\tilde{\mathcal{U}}) \in T_{\mathcal{U}}G(n,\p)$}
  \end{algorithmic}
\end{algorithm}
\tc{
 The composition $\Exp_{[U]}^{Gr} \circ \Log_{[U]}^{Gr}$ is the identity on $Gr(n,\p)$, wherever it is defined. Yet, on the level of the actual matrix representatives, the operation
 \[
  (\Exp_{[U]}^{Gr} \circ \Log_{[U]}^{Gr})([\tilde{U}_{in}]) = [\tilde{U}_{out}]
 \]
 produces a matrix $\tilde{U}_{out} \neq \tilde{U}_{in}$.
Directly recovering the input matrix can be achieved via 
a Procrustes-type preprocessing step, where $\tilde{U}$ is replaced with
$\tilde{U}_*:= \tilde{U}\Phi$, $\Phi = \argmin_{\Phi\in O(\p)}\|U-\tilde{U}\Phi\|$.
This leads to:
}
\begin{algorithm}
  \caption{Grassmann Logarithm: modified version.} 
  \label{alg:GlogMod}
  \begin{algorithmic}[1]
  \REQUIRE{base point $\mathcal{U}=[U]\in G(n,\p)$ with $U\in St(n,\p)$, 
           $\tilde{\mathcal{U}}= [\tilde{U}]\in G(n,\p)$ with $\tilde{U}\in St(n,\p)$.}
  \STATE{ $\Psi SR^T \stackrel{\text{SVD}}{:=} \tilde{U}^TU$}
  \STATE{ $\tilde{U}_* := \tilde{U}(\Psi R^T)$}
  \hfill \COMMENT{`Transition to Procrustes representative'}
  \STATE{ $L := (I -UU^T) \tilde{U}_*$ }
  \STATE{ $Q\Sigma V^T \stackrel{\text{SVD}}{:=} L$}
  \hfill \COMMENT{(thin) SVD}
  \STATE{ $\Delta := Q\arcsin(\Sigma)V^T$}
  \hfill \COMMENT{$\arcsin$ acts only on diagonal entries.}
  \ENSURE{$\Delta = \Log_\mathcal{U}^{Gr}(\tilde{\mathcal{U}}) \in T_{\mathcal{U}}G(n,\p)$}
  \end{algorithmic}
\end{algorithm}
\tc{
An additional advantage of the modified Grassmann logarithm
is that the matrix inversion $M^{-1} = (U^T\tilde{U})^{-1}$ is avoided. In fact, it is replaced by the SVD $\Psi SR^T=\tilde{U}^TU$ that is used to solve the Procrustes problem $\min_{\Phi\in O(\p)}\|U-\tilde{U}\Phi\|$.
The SVD exists also if $U^T\tilde{U}$ does not have full rank.
}
\paragraph*{Distances between subspaces}
The Riemannian logarithm provides the distance between two subspaces $\mathcal{U}=[U],\tilde{\mathcal{U}}=[\tilde{U}]\in Gr(n,\p)$
as follows: First, compute $\Delta = \Log_\mathcal{U}^{Gr}(\tilde{\mathcal{U}})$, then compute 
$\|\Delta\|_{\mathcal{U}} = \dist_{Gr}(\mathcal{U}, \tilde{\mathcal{U}})$. In practice, however, this boils down to computing the singular values of the matrix 
$M= U^T\tilde{U}$, which can be seen as follows.
By Algorithm \ref{alg:GlogMod}, 
$\|\Delta\|_{\mathcal{U}}^2 = \tr(\Delta^T\Delta) = \sum_{k=1}^p \arcsin(\sigma_k)^2$,
where the $\sigma_k$'s are the singular values of $L=(I -UU^T) \tilde{U}_*$.
These match precisely the square roots of the eigenvalues of $L^TL$.
Using the SVD of the square matrix $\tilde{U}^TU = \Psi S R^T$ as in steps 1\&2 of Algorithm \ref{alg:GlogMod}, 
the eigenvalues of $L^TL$ can be read off from
$$
  L^TL =\tilde{U}^T_*  (I -UU^T) \tilde{U}_* = I - RS^2R^T = R(I- S^{2}) R^T,
$$
so that $\sigma_k^2 = 1 - s_k^2$, when consistently ordered.
As a consequence, $s_k = \sqrt{1 -  \sigma_k^2}= \cos(\arcsin(\sigma_k))$, which implies 
\begin{equation}
\label{eq:subspaceDist}
 \dist_{Gr}(\mathcal{U}, \tilde{\mathcal{U}}) = \left(\sum_{k=1}^p \arcsin(\sigma_k)^2\right)^{\frac12}
 = \left(\sum_{k=1}^p \arccos(s_k)^2\right)^{\frac12},
\end{equation}
where $\sigma_1,\ldots,\sigma_\p$ and $s_1,\ldots,s_\p$ are the singular values of
$L$ and $\tilde{U}^TU$, respectively.

The numerical linear algebra literature knows a variety of distance measures for subspaces.
Essentially, all of them are based on the principal angles \cite[\S 2.5.1, \S 4.3]{EdelmanAriasSmith1999}.
The {\em principal angles} (or canonical angles) $\theta_1,\ldots,\theta_\p\in[0,\frac{\pi}{2}]$ between two subspaces $[U], [\tilde{U}]\in Gr(n,\p)$
are defined recursively by
\[
\cos(\theta_k) := u_k^Tv_k := 
\max_{\begin{array}{l} u\in [U], \|u\|=1\\
       u\bot u_1,\ldots, u_{k-1}
      \end{array}
     }
  \max _{\begin{array}{l} v\in [\tilde{U}], \|v\|=1\\
          v\bot v_1,\ldots, v_{k-1}
         \end{array}
        }
 u^Tv.
\]
The principal angles can be computed via
$\theta_k := \arccos(s_k) \in [0,\frac{\pi}{2}]$,
where $s_k$ is the $k$th singular value of $U^T\tilde{U}\in\R^{\p\times \p}$
\cite[\S 6.4.3]{GolubVanLoan4th}.
Hence, the Riemannian subspace distance \eqref{eq:subspaceDist} expressed in terms of the principal angles is
precisely
\begin{equation}
 \label{eq:SubspaceDist}
 \mbox{dist}([U], [\tilde{U}]) := \|\Theta\|_2, \quad \Theta = (\theta_1,\ldots,\theta_\p)\in\R^\p.
\end{equation}
In particular, \eqref{eq:SubspaceDist} shows that any two points on $Gr(n,\p)$
can be connected by a geodesic of length at most $\frac{\sqrt{\p}}{2}\pi$, see also \cite[Thm 8(b)]{Wong1967}.
%
%
%
%
\appendix
\section{Appendix}
\label{app:pre}
\paragraph*{The matrix exponential and logarithm}
The standard matrix exponential and matrix logarithm are defined via the power series
\begin{equation}
\label{eq:matrix_exp_log}
 \exp_m(X):=\sum_{j=0}^\infty{\frac{X^j}{j!}}, \quad \log_m(X):=\sum_{j=1}^\infty{(-1)^{j+1}\frac{(X-I)^j}{j}}.
\end{equation}
For $X\in\R^{n\times n}$, $\exp_m(X)$ is invertible with inverse $\exp_m(-X)$.
The following restrictions of the exponential map are important:
\[
 \exp_m|_{\sym(n)}: \sym(n)\rightarrow  SPD(n), \quad \exp_m|_{\Skew(n)}: \Skew(n)\rightarrow  SO(n).
\]
The former is a diffeomorphism \cite[Thm. 2.8]{Pennec2006}, the latter is a differentiable, surjective map 
\cite[\S. 3.11, Thm. 9]{godement2017introduction}.
For additional properties and efficient methods for numerical computation, see \cite[\S 10, 11]{Higham:2008:FM}.

A few properties of the exponential function for real or complex numbers carry over to the matrix exponential.
However, since matrices do not commute, the standard exponential law is replaced by
\begin{eqnarray}
\label{eq:BCH}
 \exp_m(Z(X,Y)) &=& \exp_m(X) \exp_m(Y),\\
 \nonumber
 Z(X,Y) &=& X + Y + \frac12 [X,Y] +\\
 \nonumber
 && \frac{1}{12}([X,[X,Y]] + [Y,[Y,X]])) - \frac{1}{24}[Y,[X,[X,Y]]] ...,
\end{eqnarray}
where $[X,Y] = XY-YX$ is the commutator bracket, or Lie bracket.
This is Dynkin's formula for the Baker-Campbell-Hausdorff series, see
\cite[\S 1.3, p. 22]{rossmann2006lie}. From a theoretical point of view,
it is important that all terms in this series can be expressed in terms of the Lie bracket.
A special case is
\[
 \exp_m(X+Y) = \exp_m(X)\exp_m(Y),\quad \mbox{if } [X,Y] = 0. 
\]

\paragraph*{Matrix square roots and the polar decomposition}
Every $S\in SPD(n)$ has a unique matrix square root in $SPD(n)$, i.e., a matrix denoted by $S^{\frac12}$ with the property
$S^{\frac12}S^{\frac12} = S$. This square root can be obtained via an EVD
$S = Q\Lambda Q^T$ by setting
\[
 S^{\frac12} := Q\sqrt{\Lambda} Q^T, 
\]
where $Q\in O(n)$, $\Lambda = \diag(\lambda_1,\ldots,\lambda_n)$ and $\lambda_i>0$ are the eigenvalues of $S$.
Every $A\in GL(n)$ can be uniquely decomposed into an orthogonal matrix times a symmetric positive definite matrix,
\[
 A = QP = Q\exp_m(X), \quad Q\in O(n), P\in SPD(n), X\in \sym(n).
\]
The polar factors can be constructed via taking the square root of the assuredly positive definite matrix $A^TA$ and subsequently setting 
$P:= (A^TA)^{\frac12}$ and $Q := AP^{-1}$.
Because the restriction of $\exp_m$ to the symmetric matrices is a diffeomorphism onto $SPD(n)$,
there is a unique $X\in\sym(n)$ with $P = \exp_m(X)$.
For details, see \cite[Thm. 2.18]{Hall_Lie2015}.
\paragraph*{The Procrustes problem}
\tc{
Let  $A,B\in \R^{n\times \p}$. The Procrustes problem aims at finding an orthogonal transformation $R^*\in O(\p)$ such that $R^*$ is the minimizer of
\[
 \min_{R\in O(\p)} \|A -BR\|_F.
\]
The optimal $R^*$ is $R^* = UV^T$, where $B^TA \stackrel{\text{SVD}}{=} U\Sigma V^T\in \R^{\p\times \p}$,
see \cite{GolubVanLoan4th}.}
%
%

\end{document}